%
\catcode`@=11
%
%
\def\bibn@me{R\'ef\'erences}
\def\bibliographym@rk{\centerline{{\sc\bibn@me}}
	\sectionmark\section{\ignorespaces}{\unskip\bibn@me}
	\bigbreak\bgroup
	\ifx\ninepoint\undefined\relax\else\ninepoint\fi}
%
%
%
\let\refsp@ce=\ 
\let\bibleftm@rk=[
\let\bibrightm@rk=]
%
%
%
\def\numero{n\raise.82ex\hbox{$\fam0\scriptscriptstyle o$}~\ignorespaces}
%
%
\newcount\equationc@unt
\newcount\bibc@unt
\newif\ifref@changes\ref@changesfalse
\newif\ifpageref@changes\ref@changesfalse
\newif\ifbib@changes\bib@changesfalse
\newif\ifref@undefined\ref@undefinedfalse
\newif\ifpageref@undefined\ref@undefinedfalse
\newif\ifbib@undefined\bib@undefinedfalse
\newwrite\@auxout
%
%
\def\eqnum{\global\advance\equationc@unt by 1%
\edef\lastref{\number\equationc@unt}%
\eqno{(\lastref)}}
%
%
%
%
%
%
\def\re@dreferences#1#2{{%
	\re@dreferenceslist{#1}#2,\undefined\@@}}
\def\re@dreferenceslist#1#2,#3\@@{\def\next{#2}%
	\expandafter\ifx\csname#1@@\meaning\next\endcsname\relax
	??\immediate\write16
	{Warning, #1-reference "\next" on page \the\pageno\space
	is undefined.}%
	\global\csname#1@undefinedtrue\endcsname
	\else\csname#1@@\meaning\next\endcsname\fi
	\ifx#3\undefined\relax
	\else,\refsp@ce\re@dreferenceslist{#1}#3\@@\fi}
%
%
%
\def\newlabel#1#2{{\def\next{#1}\newl@bel#2}}
\def\newl@bel#1#2{%
	\expandafter\xdef\csname ref@@\meaning\next\endcsname{#1}%
	\expandafter\xdef\csname pageref@@\meaning\next\endcsname{#2}}
\def\label#1{{%
	\toks0={#1}\message{ref(\lastref) \the\toks0,}%
	\ignorespaces\immediate\write\@auxout%
	{\noexpand\newlabel{\the\toks0}{{\lastref}{\the\pageno}}}%
	\def\next{#1}%
	\expandafter\ifx\csname ref@@\meaning\next\endcsname\lastref%
	\else\global\ref@changestrue\fi%
	\newlabel{#1}{{\lastref}{\the\pageno}}}}
\def\ref#1{\re@dreferences{ref}{#1}}
\def\pageref#1{\re@dreferences{pageref}{#1}}
%
%
\def\bibcite#1#2{{\def\next{#1}%
	\expandafter\xdef\csname bib@@\meaning\next\endcsname{#2}}}
\def\cite#1{\bibleftm@rk\re@dreferences{bib}{#1}\bibrightm@rk}
%
%
\def\beginthebibliography#1{\bibliographym@rk
	\setbox0\hbox{\bibleftm@rk#1\bibrightm@rk\enspace}
	\parindent=\wd0
	\global\bibc@unt=0
	\def\bibitem##1{\global\advance\bibc@unt by 1
		\edef\lastref{\number\bibc@unt}
		{\toks0={##1}
		\message{bib[\lastref] \the\toks0,}%
		\immediate\write\@auxout
		{\noexpand\bibcite{\the\toks0}{\lastref}}}
		\def\next{##1}%
		\expandafter\ifx
		\csname bib@@\meaning\next\endcsname\lastref
		\else\global\bib@changestrue\fi%
		\bibcite{##1}{\lastref}
		\medbreak
		\item{\hfill\bibleftm@rk\lastref\bibrightm@rk}%
		}
	}
\def\endthebibliography{\egroup\par}
%
%
%
\def\@closeaux{\closeout\@auxout
	\ifref@changes\immediate\write16
	{Warning, changes in references.}\fi
	\ifpageref@changes\immediate\write16
	{Warning, changes in page references.}\fi
	\ifbib@changes\immediate\write16
	{Warning, changes in bibliography.}\fi
	\ifref@undefined\immediate\write16
	{Warning, references undefined.}\fi
	\ifpageref@undefined\immediate\write16
	{Warning, page references undefined.}\fi
	\ifbib@undefined\immediate\write16
	{Warning, citations undefined.}\fi}
%
%
\immediate\openin\@auxout=\jobname.aux
\ifeof\@auxout \immediate\write16
  {Creating file \jobname.aux}
\immediate\closein\@auxout
\immediate\openout\@auxout=\jobname.aux
\immediate\write\@auxout {\relax}%
\immediate\closeout\@auxout
\else\immediate\closein\@auxout\fi
%
%
\input\jobname.aux
\immediate\openout\@auxout=\jobname.aux
%
%

\def\bibn@me{R\'ef\'erences bibliographiques}
%
\def\bibliographym@rk{\bgroup}
%
%
\outer\def\bye{ 	\par\vfill\supereject\end}

\def\Z{{\bf Z}}

\overfullrule=0pt
\def\x{{\bf  s}}

\def\bc{{\bf c}}

\def\lf{\lfloor}
\def\rf{\rfloor}
\def\lc{\lceil}
\def\rc{\rceil}
\magnification=1200

  \def\pro{\noindent {\bf{Proof : }}}

\def\house#1{\setbox1=\hbox{$\,#1\,$}%
\dimen1=\ht1 \advance\dimen1 by 2pt \dimen2=\dp1 \advance\dimen2 by 2pt
\setbox1=\hbox{\vrule height\dimen1 depth\dimen2\box1\vrule}%
\setbox1=\vbox{\hrule\box1}%
\advance\dimen1 by .4pt \ht1=\dimen1
\advance\dimen2 by .4pt \dp1=\dimen2 \box1\relax}

  \def\eps{{\varepsilon}}

\def\eps{{\varepsilon}}

\def\build#1_#2^#3{\mathrel{\mathop{\kern 0pt#1}\limitr_{#2}^{#3}}}

\def\date {le\ {\the\day}\ \ifcase\month\or 
janvier\or fevrier\or mars\or avril\or mai\or juin\or juillet\or
ao\^ut\or septembre\or octobre\or novembre\or 
d\'ecembre\fi\ {\oldstyle\the\year}}

\font\fivegoth=eufm5 \font\sevengoth=eufm7 \font\tengoth=eufm10

\newfam\gothfam \scriptscriptfont\gothfam=\fivegoth
\textfont\gothfam=\tengoth \scriptfont\gothfam=\sevengoth

\def\cqfd{\unskip\kern 6pt\penalty 500 \raise 0pt\hbox{\vrule\vbox 
to6pt{\hrule width 6pt \vfill\hrule}\vrule}\par}

\def\pro{\noindent {\it Proof. }}

\def\smallsquare{\vbox{\hrule\hbox{\vrule height 1 ex\kern 1 ex\vrule}\hrule}}
\def\cqfd{\hfill \smallsquare\vskip 3mm}

\def\bx{{\bf s}}

\def\ta{{\tilde a}}

\def\bfc{{\bf c}}
\def\bfs{{\bf s}}

\def\cK{{\cal K}}

\def\ib{b^*}   


\vskip 2mm

\centerline{\bf Combinatorial structure of Sturmian words}

\smallskip

\centerline{\bf  and continued fraction expansions of Sturmian numbers}

\vskip 13mm

\centerline{Y{\sevenrm ANN} B{\sevenrm UGEAUD} 
and M{\sevenrm ICHEL} L{\sevenrm AURENT} \footnote{}{\rm
Combinatorics on words, Sturmian word, Continued fraction expansion, 
Ostrowski numeration, Irrationality exponent.}}

{\narrower\narrower
\vskip 15mm

\proclaim Abstract. {
Let $\theta = [0; a_1, a_2, \dots]$  
be the continued fraction expansion of an irrational real number  $\theta \in (0, 1)$. 
It is well-known that the characteristic Sturmian word of slope $\theta$ is the limit 
of a sequence of finite words $(M_k)_{k \ge 0}$, with $M_k$ of length $q_k$ (the denominator of the $k$-th 
convergent to $\theta$) being a suitable concatenation of $a_k$  
copies of $M_{k-1}$ and one copy of $M_{k-2}$. 
Our first result extends this to any Sturmian word. 
Let $b \ge 2$ be an integer. 
Our second result gives the continued fraction expansion of any real number $\xi$ whose $b$-ary 
expansion is a Sturmian word $\x$ over the alphabet $\{0, b-1\}$. This 
extends a classical result of B\"ohmer who considered only the case where $\x$ is characteristic. 
As a consequence, we obtain a formula for the irrationality exponent of $\xi$ 
in terms of the slope and the intercept of $\x$. 
}

}

\vskip 6mm

\centerline{\bf 1. Introduction}

\vskip 6mm

Sturmian words are infinite words over a two letters alphabet that have exactly $n+1$ 
factors of length $n$ for every $n \ge 1$. They are the non-ultimately periodic words which are 
closest to ultimately periodic words. 
They admit several equivalent definitions and appear in many different areas of mathematics, 
including combinatorics, number theory, and dynamical systems; 
good references include Chapter 2 of \cite{Loth02}, \cite{AlSh03}, and \cite{Ar02}. 
The arithmetic description of Sturmian words is as follows. 
Throughout this paper, we let $\lfloor x \rfloor$ (resp., $\lceil x \rceil$) denote the largest (resp., smallest) 
integer less than or equal (resp., greater than or equal) to the real number $x$. 

Let $\theta$ and $\rho$ be real numbers with
$0 \le \theta, \rho < 1$ and $\theta$ irrational.
For $n \ge 1$, set
$$
s_n := s_n (\theta, \rho) = \big\lfloor n \theta + \rho \big\rfloor -
\big\lfloor (n-1) \theta + \rho \big\rfloor,
\quad
s'_n := s'_n (\theta, \rho) = \big\lceil n \theta + \rho \big\rceil -
\big\lceil (n-1) \theta + \rho \big\rceil. 
$$
Then, the infinite words
$$
{\bfs}_{\theta, \rho} := s_1 s_2 s_3 \ldots,
\quad
{\bfs}'_{\theta, \rho} := s'_1 s'_2 s'_3 \ldots 
$$
are, respectively, the lower and upper Sturmian words of slope $\theta$
and intercept $\rho$, written over the alphabet $\{0, 1\}$.
Observe that ${\bfs}_{\theta, 0}$ and ${\bfs}'_{\theta, 0}$
differ only by their first letter, thus, there exists an infinite word ${\bfc}_{\theta}$, called the 
characteristic Sturmian word of slope $\theta$,  such that 
$$
{\bfs}_{\theta, 0} = 0 {\bfc}_{\theta}, \quad  {\bfs}'_{\theta, 0} = 1 {\bfc}_{\theta}.
$$
Explicitly, we have
$$
\bfc_\theta = \bfs_{\theta,\theta}= {\bfs}'_{\theta,\theta}= c_1c_2c_3\dots ,
$$
with
$$
c_n =\lf (n+1)\theta \rf - \lf n \theta \rf = \lc (n+1)\theta \rc - \lc n \theta \rc , \quad \hbox{for $n\ge 1$.}   
$$


Alternatively, the characteristic word ${\bfc}_{\theta} = {\bfs}_{\theta, \theta} = {\bfs}'_{\theta, \theta}$ 
can be defined as follows. 
Let $[0; a_1, a_2, \ldots ]$ denote the continued fraction expansion of the slope $\theta$, 
with partial quotients $a_1, a_2, \dots$ and convergents $p_k/q_k = [0; a_1, \ldots , a_k]$ for $k \ge 1$. 
Let $(M_k)_{k \ge 0}$ be the sequence of finite words 
over the alphabet $\{{\bf a}, {\bf b}\}$ associated with $(a_j)_{j \ge 1}$ defined by
$$
M_0 = {\bf a},  \quad M_1 =  {\bf a}^{a_1 - 1} {\bf b}, \quad
M_{k} = (M_{k-1})^{a_{k}} M_{k-2}, \quad \hbox{for $k \ge 2$}. 
$$
Then, the limit $\lim_{k \to + \infty} M_k$ exists: it is the characteristic Sturmian word of slope $\theta$ over 
$\{{\bf a}, {\bf b}\}$.
Replacing ${\bf a}$ by $0$ and ${\bf b}$ by $1$, we get
$$
{\bfc}_{\theta} = \lim_{k \to + \infty} \, M_k.     \eqno (1.1)
$$
Furthermore, the length (that is, the number of letters)   
of $M_k$ is equal to $q_k$ for $k \ge 1$.

Our first result, stated as Theorem 2.1,  
extends (1.1) by showing how an arbitrary Sturmian word of slope $\theta$ and intercept $\rho$ 
can be expressed as the limit of a sequence of finite words $(V_k)_{k \ge 0}$, with $V_k$ (of length $q_k$) 
being a suitable concatenation of $a_k$ copies of $V_{k-1}$ and one copy of $V_{k-2}$, 
defined in terms of the $\theta$-Ostrowski expansion of the intercept $\rho$.

Then, we will consider some Diophantine properties of the real numbers whose sequence of 
digits in some given integer base $b$ form a Sturmian word. Such real numbers are called $b$-Sturmian
numbers, or shortly Sturmian numbers, when we do not need to refer to the base. 
The transcendence of characteristic Sturmian numbers
was established by B\"ohmer \cite{Bohm27} in 1927, 
assuming that the sequence 
of partial quotients $(a_k)_{k\ge 1}$ is unbounded. He also gave explicitly their continued 
fraction expansion; see Theorem 2.2 below. 
This has been rediscovered by Danilov \cite{Dan72}, Davison \cite{Davi77}, 
and by Adams and Davison \cite{AdDa77} 
(see also \cite{AdAl07},  Theorem 7.22 in \cite{BuLiv2}, and Section 9.3 of 
\cite{AlSh03} for a special case). 
Ferenczi and Mauduit \cite{FeMa97} used combinatorial properties of 
Sturmian words and a deep result from Diophantine approximation (Ridout's theorem, which is 
a $p$-adic extension of Roth's theorem) to establish that Sturmian numbers are transcendental. 
Specifically, they proved that every Sturmian word
contains, for some positive $\eps$, infinitely many
$(2 + \eps)$-powers of blocks (that is, a block followed by itself and by a
prefix of it   
of relative length at least $\eps$) occurring not too  
far from its beginning. 

Subsequently, Berth\'e, Holton and Zamboni \cite{BHZ06} established
that any Sturmian word, whose slope has a bounded continued fraction expansion,
has infinitely many prefixes which are $(2 + \eps)$-powers of blocks, 
for some positive real number
$\eps$ depending only on the word. This implies that the associated Sturmian 
number $\xi$ is rather close to rational numbers whose $b$-ary expansion is purely periodic 
and gives that the irrationality exponent of $\xi$ is at least equal to $2 +  \eps$.

\proclaim 
Definition 1.1.  
The irrationality exponent $\mu(\zeta)$ of an irrational real number $\zeta$ is the supremum
of the real numbers $\mu$ such that the inequality
$$
\biggl| \zeta - {p \over q} \biggr| < {1 \over q^{\mu}}
$$
has infinitely many solutions in rational numbers ${p \over q}$.
If $\mu(\zeta)$ is infinite, then $\zeta$ is called a Liouville number.

Recall that the irrationality exponent of an irrational number $\zeta$ is always 
at least equal to $2$, with equality for almost all $\zeta$, in the sense of the  
Lebesgue measure. 
As observed in \cite{Ad10} (see also Section~8.5 of \cite{BuLiv2}), 
it follows from the results of \cite{BHZ06} and \cite{AdBu11} that
the irrationality exponent of any Sturmian number exceeds $2$. 
Further progress has been made recently in \cite{BuKim19}, where 
it is proved that the irrationality exponent of a $b$-Sturmian number 
can be read on its $b$-ary expansion. 
This is equivalent to say that, among the very good 
rational approximants to a $b$-Sturmian number, 
infinitely many of them can be constructed 
by cutting its $b$-ary expansion and completing by 
periodicity. 

Furthermore, Theorem 4.3 of \cite{BuKim19} asserts that the irrationality exponent 
of a Sturmian number is at least equal to 
${5 \over 3} +  {4\sqrt{10} \over 15} = 2.5099 \ldots$, and that 
equality occurs in some cases. 
This result is obtained by means of a careful analysis of the repetitions occurring 
near the beginning of a given Sturmian word. 

Our second  main result, stated 
as Theorem 2.3, extends B\"ohmer's result and  
gives explicitly the continued fraction expansion 
of any $b$-Sturmian number over the alphabet $\{0, b-1\}$. 
From this we deduce in Theorem 2.4 an exact formula giving its irrationality exponent.  
Our approach also allows us to improve the best known transcendence measures 
for Sturmian numbers, see Theorem 2.7.

\vskip 6mm

\centerline{\bf 2. Results}

\vskip 6mm

Before stating our first result, we briefly recall the definition of 
the Ostrowski numeration system; see e.g. Proposition 2 of \cite{Ber01}. 
We keep the notation from Section 1. 
Set $q_0 = 1$ and  
$\theta_k = q_k \theta - p_k$ for $k \ge 0$. 
Note that $\theta_k < 0$ if and only if $k$ is odd. 
Let  $\sigma$ be an arbitrary number in the interval  
$[- \theta , 1 - \theta]$. Then $\sigma$   
can be  written as 
$$
\sigma = \sum_{k \ge 1} b_k \theta_{k-1},
$$
where $0 \le b_1 \le a_1 - 1$, $0 \le b_k \le a_k$ for $k \ge 2$, and 
$b_k = 0$ if $b_{k+1} = a_{k+1}$ (these are the so-called Ostrowski numeration rules).
Assume  that $\sigma$ does not belong to $\Z \theta +\Z$, 
or that $\sigma$ belongs to $\Z_{\ge 0}  \theta+\Z$. Then, we can  moreover ensure that 
 there are infinitely many odd (resp., even) 
integers $k$ such that $b_k < a_k$. 
The latter condition guarantees the unicity of the representation which is called the 
{\it Ostrowski expansion of $\sigma$}. When $\sigma$ belongs to $\Z_{\ge 0} \theta+\Z$, 
the digits $b_k$ vanish for large $k$.

\proclaim Theorem 2.1. 
Let $\theta$ and $\rho$ be real numbers with
$0 \le \theta, \rho < 1$ and $\theta$ irrational. 
Assume that $\rho$ does not belong to $\Z \theta + \Z$, or that $\rho$ 
belongs to $\Z_{\ge 1}\theta+\Z$. Then ${\bfs}_{\theta,\rho} = {\bfs}'_{\theta,\rho}$. Let 
$$
\rho-\theta = \sum_{h\ge 1} b_h \theta_{h-1}
$$
be the Ostrowski expansion of $\rho-\theta$ in base $\theta$. 
Define the words $V_{-1}, V_0, V_1, \ldots$   
by $V_{-1} = 1$, $V_0 = 0$, $V_1 = 0^{a_1-b_1-1}10^{b_1}$, and
$$
V_{k+1}= V_k^{a_{k+1}-b_{k+1}}V_{k-1}V_k^{b_{k+1}}, \quad  k \ge 1. 
$$
Then, the sequence $(V_k)_{k \ge 0}$ converges and 
$$
{\bfs}_{\theta, \rho} ={\bfs}'_{\theta,\rho}= \lim_{k \to + \infty} \, V_k.      
$$
Furthermore, setting 
$$
t_{k}= b_1+ b_2q_1+ \cdots + b_{k}q_{k-1}
\quad  \hbox{and}  \quad r_k = q_k - t_k, 
$$
and denoting by $T_k$ (resp., $R_k$) the prefix (resp., suffix) of length $t_k$ (resp., $r_k$) 
of $M_k$ for $k \ge 1$, we have 
$$
V_k = R_k T_k
\quad  \hbox{and}  \quad
M_k = T_k R_k, \quad k \ge 1.
$$

A similar result holds in the remaining case where $\rho-\theta= -m\theta +p$ 
for integers  $m\ge 1$ and $p$. This case corresponds to the sequences which are 
ultimately equal to the characteristic word ${\bfc}_\theta$. Some technical difficulties occur, 
due to the fact that the choice of the lower / upper integral part does matter; see Section 3 for a 
precise statement and its proof.

Theorem 2.1 is a key tool for our extension of the following result of B\"ohmer \cite{Bohm27}.

\proclaim Theorem 2.2 (B\"ohmer). 
For a positive real irrational number 
$\theta = [0; a_1, a_2, \ldots]$ in $(0, 1)$ and
an integer $b \ge 2$, set 
$$
\xi_b (\theta) = (b-1) \, \sum_{j= 1}^{+ \infty} \, 
{1 \over b^{\lfloor j / \theta \rfloor}}.
$$
For $k \ge 1$, let $p_k/q_k$ denote the $k$-th convergent  
to $\theta$ and set
$$
A_k := {b^{q_k} - b^{q_{k-2}} \over b^{q_{k-1}} - 1},
$$
where $q_{-1} = 0$ and $q_0 = 1$.
Then, we have
$$
\xi_b (\theta) = [0; A_1, A_2, A_3, \ldots]    
$$
and the irrationality exponent of $\xi_b (\theta)$ is given by
$$
\mu (\xi_b (\theta)) 
=  1 + \limsup_{k \to + \infty} \, {q_k  \over q_{k-1}}. 
$$

Note that $A_k$ is an integer multiple of $b^{q_{k-2}}$ since 
$q_k - q_{k-2}$ is an integer multiple of $q_{k-1}$. 

The last assertion of the theorem follows from the well-known fact that 
the irrationality exponent of an irrational real number $\zeta = [A_0; A_1, A_2, \ldots]$ is given by 
$$
\mu (\zeta) = 1 + \limsup_{j \to + \infty} \, {\log Q_{j+1} \over \log Q_j},
$$
where $[A_0; A_1, A_2, \ldots, A_j] = P_j/Q_j$, for $j \ge 1$. Indeed, the sequence 
$(P_j/Q_j)_{j \ge 1}$ comprises all the best rational approximations 
to $\zeta$ and we have 
$$
{1 \over 2 Q_{j+1} Q_j} < \biggl| \zeta - {P_j \over Q_j} \biggr| < {1 \over Q_{j+1} Q_j}.
$$

Theorem 2.2 describes the first known class of real numbers having the 
property that both their $b$-ary 
expansion (for some integer $b \ge 2$)
and their continued fraction expansion
are explicitly determined. There are only few such classes; see Section 7.6 of \cite{BuLiv2} 
for other examples. 

Our second main result extends B\"ohmer's theorem 
to an arbitrary $b$-Sturmian number with digits in $\{0, b-1\}$.
Define
$$
\xi_b (\theta, \rho) = (b-1) \, \sum_{n = 1}^{+ \infty} \, {s_n (\theta, \rho) \over b^n}, 
\quad
\xi'_b (\theta, \rho) = (b-1) \, \sum_{n = 1}^{+ \infty} \, {s'_n (\theta, \rho) \over b^n}. 
$$
Let $\xi$ denote one of these numbers. 
Let $(b_k)_{k \ge 1}$ and $(t_k)_{k \ge 1}$ be the sequences of integers defined 
in Theorem 2.1 (or in Theorem 4.2 if $\rho$ is of the form $-m \theta + p$, with $m, p$ 
nonnegative integers) applied to the Sturmian sequence defining $\xi$. Put $t_0 = 0$ and $r_0 = 1$. 
For $k \ge 0$, set 
$$
c_k = b^{r_{k} + q_{k-1}} \, 
{b^{(a_{k+1} - b_{k+1} - 1 ) q_{k}} - 1 \over b^{q_{k}} - 1},
\quad
d_k = b^{t_{k}} - 1,
$$
$$
e_k = b^{r_{k}} - 1, \quad 
f_k = b^{t_{k}} \, {b^{b_{k+1} q_{k}} - 1 \over b^{q_{k}} - 1}.
$$
We point out that some elements of these four sequences may not be positive integers. 
For example, $f_k$ is equal to $0$ when $b_{k+1} = 0$ 
and $c_{k+1}$ is equal to $0$ when $a_{k+2} = b_{k+2} + 1$. 
More intriguing is the case where $a_{k+2} = b_{k+2}$. Then, we have $b_{k+1} = 0$, thus 
$r_{k} + q_{k+1} = r_{k+1} + q_{k}$ and 
$$
c_{k+1} = b^{r_{k+1} + q_{k}} \, {b^{- q_{k+1}} - 1 \over b^{q_{k+1}} - 1} =
{b^{r_{k}} - b^{r_{k} + q_{k+1}} \over b^{q_{k+1}} - 1 } = - b^{r_{k}} = - e_{k} - 1
$$
is a negative integer. 
Keeping this in mind, and with some abuse of language, the next theorem asserts that 
$$
[0 ; c_0, d_0, 1, e_0, f_0, c_1, d_1, 1, e_1, f_1, c_2, \ldots ]
$$
is an (improper) continued fraction expansion of $\xi$. 
The precise statement is as follows.

\proclaim Theorem 2.3. 
Let $\xi$ be as above and keep the notation introduced above. 
If $a_k - b_k \ge 2$ and $b_k \ge 1$ for every $k \ge 1$,      
then the continued fraction expansion of $\xi$ is given by     
$$
\xi_b (\theta, \rho) = [0 ; c_0 + 1, e_0, f_0, c_1, d_1, 1, e_1, f_1, c_2, \ldots ].     
$$
Otherwise, let $A_1, A_2, A_3, \ldots$ be the sequence of positive integers obtained from the sequence 
$c_0, d_0, 1$, $e_0, f_0, c_1, d_1, 1, e_1, f_1, c_2, \ldots$ after the application of 
the following rules:
$$
\hbox{$(i)$ For every $k$ such that $c_{k+1} < 0$, replace the
nine integers $c_k, d_k, 1$, $e_k$, $f_k, c_{k+1}, d_{k+1},$}
$$
$$
\hbox{$1, e_{k+1}$ by the
positive integer $c_k + 1 + e_{k+1}$};
$$
$$
\hbox{$(ii)$ Replace any three consecutive elements of this new sequence of the form $x, 0, y$}
$$
$$
\hbox{by the integer $x+y$.}
$$
Then, the continued fraction expansion of $\xi$ is given by
$$
\xi_b (\theta, \rho) = [0 ; A_1, A_2, A_3, \ldots ]. 
$$

Observe that the sequence $(A_j)_{j \ge 1}$ is well-defined. Indeed, $c_k$ and $c_{k+1}$ 
cannot be both negative, since we cannot have simultaneously $a_{k+2} = b_{k+2}$ 
and $a_{k+1} = b_{k+1}$.

Let us briefly show that Theorem 2.3 includes B\"ohmer's result. 
First, note that $\xi_b (\theta) = \xi_b (\theta, \theta)$, since, for a positive integer $j$, we have 
$\lfloor j / \theta \rfloor$ equals the integer $\ell$ if and only if $\ell < j / \theta < \ell + 1$, that is, 
if and only if, $\lfloor (\ell + 1)  \theta \rfloor - \lfloor \ell \theta \rfloor = 1$. Then, observe that 
the Ostrowski expansion of $\theta - \theta = 0$ in base $\theta$ is given by the constant
sequence equal to $0$. Consequently, the sequences defined in Theorem 2.3 are 
equal to
$$
d_k = f_k = 0, \ \ e_k = b^{q_k} - 1, \ \ 
c_k = b^{q_k + q_{k-1}} \, {b^{(a_{k+1}  - 1 ) q_{k}} - 1 \over b^{q_{k}} - 1}, \quad  k \ge 0.    
$$
It then follows from Theorem 2.3 that 
$$
\eqalign{
\xi_b (\theta) 
& = \biggr[0; b^{q_0 + 0} \, {b^{(a_{1}  - 1 ) q_{0}} - 1 \over b^{q_{0}} - 1}, 0, 1, b^{q_0} - 1, 0, 
b^{q_1 + q_{0}} \, {b^{(a_{2}  - 1 ) q_{1}} - 1 \over b^{q_{1}} - 1}, 0, 1, b^{q_1} - 1, 0, c_2, \ldots \biggr] \cr
& = \biggr[0; b^{q_0 + 0} \, {b^{(a_{1}  - 1 ) q_{0}} - 1 \over b^{q_{0}} - 1} + 1, b^{q_0} - 1 + 
b^{q_1 + q_{0}} \, {b^{(a_{2}  - 1 ) q_{1}} - 1 \over b^{q_{1}} - 1} + 1, b^{q_1} - 1, 0, c_2, \ldots \biggr]  \cr
& = \biggr[0;   {b^{a_{1} q_{0}} - 1 \over b^{q_{0}} - 1},  
{b^{q_{2}} - b^{q_{0}} \over b^{q_{1}} - 1}, b^{q_1} - 1, 0, c_2, \ldots \biggr] \cr
& = \biggr[0;   {b^{q_{1}} - 1 \over b^{q_{0}} - 1},  
{b^{q_{2}} - b^{q_{0}} \over b^{q_{1}} - 1}, {b^{q_{3}} - b^{q_{1}} \over b^{q_{2}} - 1}, \ldots \biggr]. \cr
}
$$
We get the sequence of partial quotients $c_0+1, e_0 + c_1 + 1, e_1 + c_2 + 1, \ldots$ and we recover 
Theorem 2.2. 

Theorem 2.3 is proved in Section 7, where we give additional informations 
on the shape of the convergents to $\xi$ and its partial quotients; see Proposition 7.2.

\medskip

As a consequence of Theorem 2.3, we obtain an expression for the 
irrationality exponent of any Sturmian number in terms of its slope and its intercept. 

Keep our notation and define
$$
\nu_k(1) = 2+ {t_k\over r_{k+1}}, \quad \nu_k(2)
= 2+ {r_{k} \over r_{k+1}+ t_k}, 
$$
$$
\nu_k(3) = 1+ { q_{k+1}\over r_{k+1}+ q_k},
\quad \nu_k(4)= 1+ { r_{k+2}\over q_{k+1}}.
$$
Put 
$$
\eqalign{
\nu(1) & = \limsup_{k \to + \infty} \, \{ \nu_k(1) \, : \, a_{k+1}-b_{k+1} \ge 1 
\hbox{ and $a_{k+2}-b_{k+2} \ge 1$}\},
\cr
\nu(2)  &= \limsup_{k \to + \infty} \, \{ \nu_k(2) \, : \, a_{k+2}-b_{k+2} \ge 1 \},
\cr}
$$
and, for $j = 3, 4$, 
$$
\nu(j) = \limsup_{k \to + \infty}  \nu_k(j).
$$

\proclaim 
Theorem 2.4. 
Let $\xi$ be as above. 
Then, its irrationality exponent is equal to
$$
\max\{ \nu (1), \nu (2), \nu (3), \nu (4) \}. 
$$

We recover, for the initial repetitions, the formulas found in 
\cite{BHZ06} for the critical initial exponent, namely the contributions of $\nu (2)$ and $\nu (4)$. 
Theorem 2.4 is established at the end of Section 6; see Theorem 6.3.

Furthermore, we derive easily a necessary and sufficient condition under which a 
Sturmian number is a Liouville number, thereby reproving the first part of 
Th\'eor\`eme 3.1 of \cite{AdBu11} (see also \cite{Kom96}). 

\proclaim
Corollary 2.5. 
A Sturmian number is a Liouville number if and only if its slope has 
unbounded partial quotients in its continued fraction expansion.

Theorem 2.4 allows us to study in depth the irrationality exponents of Sturmian numbers. 
For instance, we can fix a slope $\theta$ and consider the spectrum ${\cal L} (\theta)$ 
consisting of the set the irrationality exponents of Sturmian numbers of slope $\theta$. 

\proclaim 
Theorem 2.6. 
Let $\theta$ be an irrational number in $(0, 1)$ with bounded partial quotients. 
Then, 
$$
{\cal L} (\theta) \subset \biggl[{5 \over 3} +  {4\sqrt{10} \over 15}, 1 + \mu(\xi_b(\theta))\biggr]
$$
and there exists an intercept $\rho (\theta)$ such that 
$$
\mu(\xi_b(\theta, \rho (\theta))) = 1 + \mu(\xi_b(\theta)). 
$$

A detailed study of the sets ${\cal L} (\theta)$ will be the purpose of a forthcoming work.

Theorem 2.3 allows us also to improve the best known 
transcendence measures for Sturmian numbers. 
Let $\zeta$ be a transcendental real number. 
Following Koksma \cite{Ko39}, for any integer $d \ge 1$, we denote by
$w_d^* (\zeta)$ the supremum of the exponents $w$ for which
$$
0 < |\zeta - \alpha| < H(\alpha)^{-w-1}
$$
has infinitely many solutions in real algebraic numbers $\alpha$ of
degree at most $d$. Here, $H(\alpha)$ stands for the 
na\"\i ve height of the minimal defining polynomial of $\alpha$
over $\Z$. Clearly, the functions $\mu - 1$ and $w_1^*$ are equal   
and the functions $w_d^*$ are invariant by rational translation and by    
multiplication by a nonzero rational number, for $d \ge 1$. 
We direct the reader to \cite{BuLiv1} for classical results on the functions $w_d^*$
and on Mahler's and Koksma's classifications of real numbers. 
As a particular case of 
Th\'eor\`eme 1.1 of \cite{AdBu11}, we know that, for any Sturmian number $\xi$ which is not 
a Liouville number, there exists a positive real number $c$, depending only on $\xi$, such that 
$$
w_d^* (\xi) \le (2d)^{c (\log 3d) (\log \log 3d)}, \quad d \ge 1.
$$
This can be improved as follows. 

\proclaim 
Theorem 2.7. 
Let $\xi$ be a Sturmian number. 
Assume that the partial quotients of its slope are ultimately bounded from above by $M$. 
Then, there exists a positive real number $\kappa$, depending only on $M$, such that 
$$
w_d^* (\xi) \le (2d)^{\kappa (\log \log 3d)}, \quad d \ge 1.
$$

We point out that the transcendence measure obtained in Theorem 2.7 does not depend 
on the intercept of the Sturmian number. 

We believe  
that Theorem 2.1 will have many applications. 
We will use it in a follow-up work devoted to the transcendence of Hecke--Mahler series 
evaluated at algebraic points. We refer
to  \cite{BKLN,LaNo18,LaNo21} for various applications of Sturmian numbers to the dynamics of piecewise affine maps.

The present paper is organized as follows. 
We show in Section 3 that any Sturmian word $\x$ of slope $\theta$ 
and intercept $\rho$ can be expressed in a similar way 
as in (1.1) and we define its formal intercept. The link between the formal intercept and the 
expansion of the intercept $\rho$ in the $\theta$-Ostrowski numeration system is 
established in Section 4, thereby proving Theorem 2.1. In Section 5, we apply Theorem 2.1 to 
give a precise description of the repetitions occurring near the beginning of $\x$. 
From this, in the next section, 
we deduce four one-parametric families of rational numbers which approximate very well
the Sturmian number $\xi$ associated to $\x$, the exact rate
of approximation to $\xi$ by these rational numbers 
being given in Theorem 6.1.  
We derive   
the continued fraction expansion of $\xi$ 
in Section 7, thereby proving Theorems 2.4 and 2.5, since we see that 
all the very good approximants to $\xi$ belong to one of the four families defined in Section 6. 
The final Section is devoted to the proofs of the other results stated in Section 2.

\vskip 5mm

\centerline{\bf 3. The formal intercept of a Sturmian word}

\vskip 5mm

We keep the notation of Section 1 with the alphabet $\{0,1\}$.   
Let $\bx$ be an arbitrary Sturmian word of slope $\theta$. 
The goal of this section is to establish that any Sturmian word can be expressed as in (1.1), that is, 
as the limit of a suitable sequence $(V_k)_{k \ge 1}$ of binary words $V_k$
of length $q_k$ constructed inductively.     

Throughout, the length $|W|$ of a finite word $W$, that is, the number of letters composing $W$, is denoted 
by $|W|$. If $W$ has at least one letter (resp., at least two letters), then $W^-$ (resp., $W^{--}$) 
denotes the work $W$ deprived of its last letter (resp., its last two letters).   

\proclaim 
Definition  3.1. 
A word  $V$ is a {\it conjugate} of $M_k$ if there exist words $T$ and $R$ such that 
$$
V = R T \quad  \hbox{and}  \quad M_k= T R,
$$
with $0 \le t : =  | T| < q_k$. 
Then, $R$ is the non-empty suffix of $M_k$ of length $q_k-t$. 

Observe that the $q_k$   conjugates $V$ of the word $M_k$
are distinct.  We label these translated words $V$  by the length $t, 0 \le t< q_k$ 
of the (possibly empty) prefix $T$ in the decomposition $M_k= TR, V=RT$. 
The whole set of  conjugates  $V$ of $M_k$ is clearly obtained as the set of factors 
of length $q_k$ in the word $M_kM_k^-$. Each such factor $V$ is determined by its $q_k-1$ first letters 
which form the $q_k$ distinct factors of length $q_k-1$ contained in the 
word $M_k M_k^{--}$. 

 As an example, for $k=1$, we have $M_1 = 0^{a_1-1}1$. Any conjugate $V$ of $M_1$
  can be written in the form
$$
V= 0^{a_1-1 -b_1}10^{b_1} = RT,  \quad M_1 = TR, \quad{\rm with } \quad T= 0^{b_1},  
\quad R= 0^{a_1-1-b_1}1, 
$$
for  some integer $b_1$ with $0\le b_1\le a_1-1$. Thus, in this case, we have $t=b_1$.

\proclaim
Definition 3.2. 
For each $k\ge 1$, let $V_k$ be the conjugate of $M_k$ 
whose first $q_k-1$ letters coincide with those of $\bx$. Let $T_k$ and $R_k$ be the words such that 
$$
V_k = R_k T_k \quad { \hbox{and} } \quad M_k = T_k R_k,
$$
with $R_k$ non-empty. Denote by $t_k$ the length of $T_k$. 
Put $R_{-1} = 1$, $R_0 = 0$, and let $T_{-1}$ and $T_0$ be the empty word.    

Then, the following recursion formulae hold. 
The notion of formal intercept was first introduced by Wojcik \cite{Wo19}, 
but our presentation is different.     

\proclaim 
Lemma 3.3 (formal intercept). 
Put $t_1 = \ib_1$.
For any $k\ge 1$, there exists an
integer $\ib_{k+1}$ such that $0\le \ib_{k+1} \le a_{k+1}$ and
$$
t_{k+1}= t_k + \ib_{k+1}q_k. 
$$
When $\ib_{k+1}= a_{k+1}$, we necessarily have $t_k < q_{k-1}$,  
so that  $\ib_k=0 $ and  $t_k=t_{k-1}$ in this case. 
Moreover, the sequences of words $(T_k)_{k \ge 0}$ and $(R_k)_{k \ge 0}$    
satisfy the recursion formulae
$$
T_{k+1}= M_k^{\ib_{k+1}}T_k = T_kV_k^{\ib_{k+1}}
$$
and
$$
R_{k+1}= \cases{ R_k M_k^{a_{k+1}-\ib_{k+1}-1}M_{k-1} & if $\ib_{k+1}< a_{k+1}$, 
\cr
R_{k-1} & if $\ib_{k+1}=a_{k+1}$,
\cr}
$$
for $k \ge 0$.  
The sequence $(\ib_k)_{k \ge 1}$ is called the formal intercept of $\bx$. 

\pro
The word $V_{k+1}$ is a factor of the word
$
M_{k+1}M_{k+1}^- 
$
beginning  somewhere on the first factor $M_{k+1}$. Assume first that $V_{k+1}$ 
begins on the prefix $M_k^{a_{k+1}}$ of $M_{k+1}= M_k^{a_{k+1}}M_{k-1}$ and let $P $  
be the prefix of length $q_k$ of $V_{k+1}$. Thus, for some integer $0 \le \ib_{k+1}< a_{k+1}$, 
the prefix $P$ begins on the second factor $M_k$ in the product 
$M_k^{a_{k+1}}= M_k^{\ib_{k+1}}M_k M_k^{a_{k+1}-\ib_{k+1}-1}$. Then,  $P$ is a factor of
$$
M_k M_k^{a_{k+1}-\ib_{k+1}-1}M_{k+1}^-= 
M_k M_k^{2 a_{k+1}-\ib_{k+1}-1}M_{k-1}^-
$$
beginning on the first factor $M_k$. Since $2 a_{k+1}-\ib_{k+1} -1\ge 1$, we see that $P$ 
is located over the product $M_kM_k^-$, where $M_k^-$ 
is the prefix of $M_k^{2 a_{k+1}-\ib_{k+1}-1}$ of length $q_k-1$. 
As the first $q_k-1$ letters of $P$ coincide with those of $\bx$, we deduce that $P= V_k=R_kT_k$,  
and next that 
$$
T_{k+1}=M_k^{\ib_{k+1}}T_k 
\quad \hbox{and} \quad
R_{k+1}=  R_k M_k^{a_{k+1}-\ib_{k+1}-1}M_{k-1}.
$$
Note finally that
$$
M_k^{\ib_{k+1}}T_k = (T_kR_k)^{\ib_{k+1}}T_k= T_k(R_kT_k)^{\ib_{k+1}}= T_kV_k^{\ib_{k+1}}.
$$

Suppose now that $V_{k+1}$ begins on the second   factor $M_{k-1}$ in 
$$
M_{k+1}M_{k+1}^-= M_k^{a_{k+1}}M_{k-1}M_{k+1}^-= M_k^{a_{k+1}}T_{k-1}R_{k-1}M_{k+1}^-
$$ 
and put $\ib_{k+1}= a_{k+1}$. Then, 
$$
T_{k+1} = M_k^{a_{k+1}}T_{k-1}
\quad{\rm and}\quad R_{k+1}= R_{k-1},
$$
observing   that $V_{k-1}=R_{k-1}T_{k-1}$ equals  the prefix of $V_{k+1}$ of  length  $q_{k-1}$.  
Notice now that $M_k$ is a prefix of  
$M_{k-1}M_{k+1}^-$. Writing 
$$
M_{k-1} M_{k+1}^{- }= M_k\cdots = T_{k-1}R_{k-1}M_{k-1}^{a_k-1}M_{k-2}\cdots
$$
 we see that $T_k= T_{k-1}$ and $R_k = R_{k-1}M_{k-1}^{a_k-1}M_{k-2}$.  
 Thus $\ib_k=0$ by the preceding case applied to  the level $k-1$.   
\cqfd

We now deal with binary recursions  expressing  $V_{k+1} $ in terms of $V_k$ and $V_{k-1}$ 
extending the classical formulae    
$ M_{k+1}= M_k^{a_{k+1}}M_{k-1}$. 
Set $V_{-1} = R_{-1} T_{-1} = 1$ and $V_0 = R_0 T_0 = 0$.     

\proclaim
Lemma 3.4 (binary recursion). We have the relation 
$
V_1 = V_0^{a_{1} -1- \ib_{1}} V_{-1} V_0^{\ib_{1}},
$
while for any $k\ge 1$, we have
$$
V_{k+1}= V_k^{a_{k+1} - \ib_{k+1}} V_{k-1} V_k^{\ib_{k+1}}.
$$

\pro 
The expression 
$$
V_1 = 0^{a_1-1-\ib_1}10^{\ib_1}
$$
yields obviously the relation for $V_1$. 

For $k\ge 1$, we distinguish two cases, either $\ib_{k+1} < a_{k+1}$ or $\ib_{k+1}= a_{k+1}$. 
Assume first that $\ib_{k+1} < a_{k+1}$. According to Lemma 3.3, 
we write $ V_{k+1} = R_{k+1} T_{k+1}$ with
$$
\eqalign{
R_{k+1} &= R_k M_k^{a_{k+1} - \ib_{k+1} - 1} M_{k-1} = 
R_k(T_kR_k)^{a_{k+1}-\ib_{k+1}-1} T_{k-1}R_{k-1} 
\cr
&= (R_kT_k)^{a_{k+1}-\ib_{k+1}-1}R_kT_{k-1}R_{k-1}= V_k^{a_{k+1}-\ib_{k+1}-1} R_kT_{k-1}R_{k-1}
}
$$
and 
$$
T_{k+1} = T_k V_k^{\ib_{k+1}}. 
$$
Thus
$$
V_{k+1} =  V_k^{a_{k+1} - \ib_{k+1} - 1} R_k T_{k-1} R_{k-1} T_k V_k^{\ib_{k+1}}.
$$
Since
$$
R_k T_{k-1} R_{k-1} T_k = R_k T_{k-1} R_{k-1} T_{k-1} (R_{k-1} T_{k-1})^{\ib_k} 
= R_k T_k R_{k-1} T_{k-1}=V_kV_{k-1},
$$
we get
$$
V_{k+1} = V_k^{a_{k+1} - \ib_{k+1}} V_{k-1} V_k^{\ib_{k+1}}.
$$

Assume now that $\ib_{k+1} = a_{k+1}$. Then $\ib_k=0$. 
From  Lemma 3.3, we know that $T_k = T_{k-1}$ and that
$$
T_{k+1}= T_kV_k^{\ib_{k+1}} = T_{k-1}V_k^{\ib_{k+1}}
\quad{\rm and }\quad
R_{k+1}= R_{k-1}.
$$
Thus 
$$
V_{k+1} = R_{k+1} T_{k+1} = R_{k-1} T_{k-1} V_k^{\ib_{k+1}} = V_{k-1} V_k^{\ib_{k+1}},
$$
as asserted.
\cqfd

\medskip

We conclude this section with a corollary, which shows how any prefix of $M_{n+1}$
can be expressed in terms of $M_0, \ldots , M_n$.  
 
Recall that the Ostrowski numeration system in base $\theta$ is defined as follows:  
every positive integer $N$ can be uniquely written in the form
$$
N = d_1 + d_2 q_1 + \ldots + d_{r+1} q_r,
$$
where $0 \le d_j \le a_j$ for $j = 1, \ldots , r+1$, $d_{r+1} > 0$, $d_1 < a_1$ and 
$d_j = 0$ if $d_{j+1} = a_{j+1}$.

\proclaim
Corollary 3.5 (product formula for prefixes). Let $T$ be the prefix of $M_{n+1}$ of length $t <q_{n+1}$.
Write 
$$
t = d_1 + d_2q_1+ \cdots + d_{n+1} q_n
$$
where $d_1, \dots , d_{n+1}$ are the digits of the integer $t$ in the Ostrowski numeration 
system in base $\theta$. 
Then, we have the product formula
$$
T = M_n^{d_{n+1}}M_{n-1}^{d_n}\cdots M_0^{d_1} = V_0^{d_1} V_1^{d_2}\cdots V_n^{d_{n+1}},
$$
where the words $V_0, \dots, V_n$ are defined recursively by the formulae
$$
 V_0 = 1,\, \,V_1 = 0^{a_1-d_1-1}10^{d_1}, \, \,  
 V_{k+1}= V_k^{a_{k+1}-d_{k+1}}V_{k-1}V_k^{d_k},\, \,  1 \le k < n.
$$

\pro
By Lemma 3.3, we have $T= T_{n+1}$ and $t=t_{n+1}$. The recurrence relations
$$
T_{k+1}= M_k^{d_{k+1}} T_k= T_k V_k^{d_{k+1}}
$$
yield  inductively the product formula
$$
T= T_{n+1} = M_n^{d_{n+1}}M_{n-1}^{d_n}\cdots M_0^{d_1} = V_0^{d_1} V_1^{d_2}\cdots V_n^{d_{n+1}}. 
$$
This establishes the corollary.     
\cqfd


\vskip 5mm

\centerline{\bf 4. Linking formal intercept and Ostrowski numeration}

\vskip 5mm

We link the formal intercept, that is the sequence $(\ib_k)_{k\ge 1}$ such that 
$$
t_{k}= \ib_1+ \ib_2q_1+ \cdots + \ib_{k}q_{k-1}, \quad   k\ge 1, 
$$
 to the intercept $\rho$ thanks to the 

\proclaim
Proposition 4.1. 
Let $0< \rho < 1$ be a real number either not belonging to $\Z \theta + \Z$,  
or of the form $\Z_{\ge 1}\theta + \Z$. Let 
$$
\rho-\theta = \sum_{h\ge 1}b_h\theta_{h-1}
$$
be the Ostrowski expansion of $\rho-\theta$ in base $\theta$. 
For every $k\ge 1$, put 
$$
t_{k}= b_1+ b_2q_1+ \cdots + b_{k}q_{k-1}.
$$
Then, $t_k$ is the length of the word $T_k$ associated to  the Sturmian word 
$ \bfs_{\theta,\rho}= \bfs'_{\theta,\rho}$. In other words, we have $b_k= \ib_k$ for $k \ge 1$, meaning 
that the formal intercept of this Sturmian word coincides with the sequence of digits 
of the number $\rho-\theta$ 
in its Ostrowski expansion in base $\theta$.

\pro
By definition, we have
$$
s_n = \lf n\theta + \rho \rf -\lf (n-1)\theta + \rho \rf, \quad  n\ge 1,
$$
and 
$$
s'_n = \lc n\theta + \rho \rc -\lc (n-1)\theta + \rho \rc, \quad  n\ge 1,
$$
while the $n$-th letter of $\bfc_\theta$ is 
$$
 c_n= \lf (n+1)\theta \rf -\lf n\theta\rf= \lc (n+1)\theta \rc -\lc n\theta\rc , \quad  n\ge 1. 
 $$
 Thus
$$
s_n = \lf (n+1)\theta + \rho-\theta \rf -\lf n\theta + \rho-\theta \rf 
= \lf (n+1+t_k )\theta +  \sigma_k \rf - \lf (n+t_k)\theta + \sigma_k \rf,
$$
and 
$$
s'_n = \lc (n+1)\theta + \rho-\theta \rc -\lc n\theta + \rho-\theta \rc 
= \lc (n+1+t_k )\theta +  \sigma_k \rc - \lc (n+t_k)\theta + \sigma_k \rc,
$$
where we have set
$$
\sigma_k= \sum_{h\ge k} b_{h+1}\theta_h.
$$
We claim that
$$
\lf q\theta + \sigma_k \rf = \lf q\theta \rf
\quad{\rm and}\quad 
\lc q\theta + \sigma_k \rc = \lc q\theta \rc
$$
for every integer $q$ with $1 \le q \le q_k + t_k$. This yields that
$$
s_n=  \lf (n+1+t_k )\theta +  \sigma_k \rf - \lf (n+t_k)\theta + \sigma_k \rf
 =
  \lf (n+1+t_k )\theta \rf - \lf (n+t_k)\theta  \rf  = c_{n+t_k}
  $$
  and 
  $$
s'_n=  \lc (n+1+t_k )\theta +  \sigma_k \rc - \lc (n+t_k)\theta + \sigma_k \rc
 =
  \lc (n+1+t_k )\theta \rc - \lc (n+t_k)\theta  \rc  = c_{n+t_k}
  $$
for every $1 \le n \le q_k-1$, and  will establish the proposition, noting that $M_k M_k$ 
is a prefix of $\bc_\theta$.

To that purpose, we bound $| \sigma_k |$. Observe that $\theta_k$ is positive when $k$ is even and negative
when $k$ is odd. Moreover $b_{h+1}\le a_{h+1}$ for any $h\ge 1$, while $b_1 \le a_1-1$.  Thus,
$$
\eqalign{
| \sigma_k | & < \max(\vert a_{k+1} \theta_k  + a_{k+3} \theta_{k+2}+ \cdots \vert, \vert a_{k+2} \theta_{k+1} 
 + a_{k+4} \theta_{k+3}+ \cdots \vert )  \cr
& = \max ( | \theta_{k-1} | , | \theta_k |) = | \theta_{k-1}|, \cr
}
$$  
 noting that
 $$
\displaylines
{
a_{k+1} \theta_k  + a_{k+3} \theta_{k+2}+ \cdots = \lim_{n\rightarrow \infty}
\biggl( (\sum_{h=0}^n a_{k+2h+1}q_{k+2h})\theta -(\sum_{h=0}^n a_{k+2h+1}p_{k+2h}) \biggr)
 \cr
 = \lim_{n\rightarrow \infty}\biggl( \sum_{h=0}^n (q_{k+2h+1}-q_{k+2h-1})\theta 
 - \sum_{h=0}^n (p_{k+2h+1}-p_{k+2h-1})\biggr) 
 \cr
 = \lim_{n\rightarrow \infty} \Big((q_{k+2n+1} -q_{k-1})\theta -(p_{k+2n+1}-p_{k-1})\Big) = -\theta_{k-1}.
 }
 $$
 The inequality $| \sigma_k | < | \theta_{k-1}| $ is strict because either $ \rho -\theta$ 
 does not belong to $\Z \theta + \Z$, or $\rho-\theta$ belong to  $\Z_{\ge 0} \theta+\Z$, 
 so that the sequence of digits $(b_h)_{h\ge 1}$ 
 cannot be ultimately of the form $ a_{k+1}, 0, a_{k+3}, 0 \dots$.
 
Observe now that 
$$
| \sigma_k - b_{k+1}\theta_k | < a_{k+2}| \theta_{k+1} | + a_{k+4}| \theta_{k+3}|  + \cdots  = | \theta_k |.
$$
 It follows that $\sigma_k$ and $\theta_k$ share the same sign when $b_{k+1}\ge 1$ 
 and that $| \sigma_k | < | \theta_k | $ when $b_{k+1}=0$. In particular,  
 the stronger inequality $| \sigma_k | < | \theta_k | $ holds when $\theta_k$ and $\sigma_k$ 
 have opposite signs.
 
 The  upper bound $| \sigma_k | < | \theta_{k-1} |$ can also be sharpened when $t_k \ge q_{k-1}$. 
 Indeed in this case we have $b_k \ge 1$ and thus $b_{k+1}$ cannot be equal to $a_{k+1}$ by 
 Ostrowski's numeration rules. We now bound $b_{k+1}\le a_{k+1}-1$ to obtain
 $$
 | \sigma_k | <   | \theta_{k-1} | - | \theta_k |.\eqno{(4.1)}
 $$
  
 Denote by $\|  x \|$ the distance from the real number $x$ to the closest integer. 
 We now show that $\| q\theta\| $ is larger than $| \sigma_k |$ when $q$ differs from $q_k$, 
 so that $q\theta $ and $q\theta +\sigma_k$ belong to the same integer open interval of length $1$ 
 and have thus the same upper and lower integer parts.  We distinguish three cases. If $q< q_{k}$, then 
 $$
 \| q\theta\|  \ge | \theta_{k-1} | > | \sigma_k |, 
 $$
 as required.  Assume secondly that $t_k < q_{k-1}$ and $q = q_k +v $ for some $1\le v \le t_k$. Then 
 $$
 \| q \theta \| = \| v \theta + \theta_k\| \ge \| v \theta \| - | \theta_k | \ge | \theta_{k-2}| - | \theta_k| 
 = | \theta_{k-2} -  \theta_k|  \ge | \theta_{k-1} | > | \sigma_k |.
$$
Thirdly, assume $t_k \ge q_{k-1}$ and $q= q_k + v$ with $1\le v \le q_k-1$. Then,
$$
\| q \theta \| \ge \| v \theta \| - | \theta_k | \ge | \theta_{k-1}| - | \theta_k|   > | \sigma_k |,
$$
by $(4.1)$. 

These three cases cover all the values of $q$ with $1 \le q \le q_k + t_k$, except $q = q_k$, which we
consider now. 
We have $\| q_k\theta \| = | \theta_k |$. 
When $\theta_k$ and $\sigma_k$ share the same sign, we have 
$$
| \theta_k | < | \theta_k +\sigma_k | = | \theta_k | + | \sigma_k | \le | \theta_k | + | \theta_{k-1}| <1.
$$
Thus, $q_k\theta$ and $q_k\theta+ \sigma_k$ both belong either to $(p_k,p_k+1)$ or to $(p_k-1,p_k)$.
When $\theta_k$ and $\sigma_k$ 
have opposite signs, we know that $| \sigma_k | < | \theta_k | $, so that $\theta_k $ and $\theta_k + \sigma_k$ 
have the same sign and both have absolute value less than $1$. The claim is proved, which yields the proposition.
\cqfd

\bigskip

A similar result holds in the remaining case where $\rho-\theta= -m\theta +p$ for integers  $m\ge 1$ and $p$.  
Assume first that $\rho$ is positive, that is to say $m\ge 2$. 
Let $l \ge 0$ be defined by the inequalities $q_l <  m  \le q_{l+1}$ and let
$$
q_{l+1} -m = b_1q_0 + \cdots + b_{l+1}q_l
$$
be the Ostrowski expansion of the integer $q_{l+1}-m$ (see the definition at the end of Section 3). 
Observe that 
$$
b_{l+1}\le a_{l+1}-1 \quad {\rm and \, \, that}\quad  b_l=0 \quad {\rm when} \quad  b_{l+1}= a_{l+1}-1.
$$
Then $ \rho-\theta\in (-\theta,1-\theta) $ has two Ostrowski expansions of the form
$$
\rho-\theta  =  b_1\theta_0 + \cdots + b_{l+1}\theta_l + \sum_{k\ge 1}a_{l+2k+1}\theta_{l+2k} 
$$
and 
$$
\rho-\theta=
b_1\theta_0 + \cdots + b_{l}\theta_{l-1} +(b_{l+1}+1)\theta_l +(a_{l+2}-1)\theta_{l+1}  + \sum_{k\ge 2}a_{l+2k}\theta_{l+2k-1},
$$
when $l\ge 1$, or 
$$
\rho-\theta=
(b_1+1)\theta_0 +(a_{2}-1)\theta_{1}  + \sum_{k\ge 2}a_{2k}\theta_{2k-1},
$$
when $l=0$. 
Set 
$$
b_{l+2}= 0, b_{l+3} = a_{l+3},b_{l+4}= 0, b_{l+5}= a_{l+5}, \dots
$$
and 
$$
b'_1= b_1,\dots , b'_l=b_l, b'_{l+1}= b_{l+1}+1, b'_{l+2}= a_{l+2}-1, b'_{l+3}=0, b'_{l+4}= a_{l+4}, \dots
$$
when $l\ge 1$, or 
$$
b'_1= b_1+1, b'_{2}= a_{2}-1, b'_{3}=0, b'_{4}= a_{4}, \dots
$$
when $l=0$, so that $(b_k)_{k\ge 1}$ and $(b'_k)_{k\ge 1} $ are the  sequences of digits appearing in the two
 above expansions of $\rho-\theta$. Notice that both sequences satisfy the Ostrowski numeration rules
  for digits in base $\theta$. When $\rho=0$, we use the two proper expansions
$$
1-\theta = (a_1-1)\theta_0 + \sum_{k\ge 1}a_{2k+1}\theta_{2k},
$$
and 
$$
-\theta = \sum_{k\ge 1}a_{2k}\theta_{2k-1},
$$
to define respectively the sequences of digits  $(b_k)_{k\ge 1}$ and $(b'_k)_{k\ge 1} $.
Then, we have the following analogue of Theorem 2.1. 

\proclaim
Theorem 4.2. 
Assume that $\rho -\theta = -m \theta +p$ where $m\ge 1 $ and $p$ are integers. 
When $m\ge 2$, let $l \ge 0$ be defined by the inequalities $q_l <  m  \le q_{l+1}$. When $m=1$, set $l=0$.   
Let $(V_k)_{k\ge 0} $ and $(V'_k)_{k\ge 0}$ be the two sequences of words recursively 
defined as in Theorem 2.1, 
with respect to the two sequences of digits $(b_k)_{k\ge 1}$ and $(b'_k)_{k\ge 1}$
defined above. 
When  $l$ is odd, we have
$$
{\bfs}_{\theta, \rho} = \lim_{k \to + \infty} \, V_k \quad {\rm and}\quad
{\bfs}'_{\theta, \rho} = \lim_{k \to + \infty} \, V'_k   .   
$$
When $l$ is even, we have
$$
{\bfs}_{\theta, \rho} = \lim_{k \to + \infty} \, V'_k \quad {\rm and}\quad
{\bfs}'_{\theta, \rho} = \lim_{k \to + \infty} \, V_k   .   
$$
Moreover, the analogous decompositions $V_k= R_kT_k$ and $V'_k= R'_kT'_k$, 
as in Theorem 2.1,  hold true with 
$$
t_k= b_1+ \cdots + b_kq_{k-1}\quad {\rm and} \quad t'_k= b'_1+ \cdots + b'_kq_{k-1}.
$$

\pro
We  only give a complete proof for the sequence of digits 
$$
(b_k)_{k\ge 1} = \Big\{ b_1, \dots , b_{l+1},
 0, a_{l+3}, 0, a_{l+5}, \dots\Big\} .
$$
Assume that $m\ge 2$ and $l$ is odd,  and recall the notations
 $$
 t_k = b_1 + \cdots + b_k q_{k-1} \quad { \rm and} \quad \sigma_k = b_{
k+1}\theta_k + b_{k+2}\theta_{k+2}+ \cdots
$$
  The argumentation is similar to the proof of Proposition 4.1. It suffices to show that
$$
\lf q\theta + \sigma_k \rf = \lf q\theta \rf  \eqno{(4.2)}
$$
for every integer $q$ with $1 \le q \le q_k + t_k$.  If $k \le l$, we compute
 $$
 \sigma_k = b_{
k+1}\theta_k + \cdots + b_{l+1}\theta_{l+1} -\theta_{l+1}.
$$
Since the tail  $b_{k+1}, b_{k+2}\dots$ of the sequence $(b_k)_{k \ge 1}$
contains the subsequence $\dots b_{l+1},0, a_{l+3}, \dots$ and that $b_{l+1}\le a_{l+1} -1$, 
observe that this tail is neither  of the form $a_j,0,a_{j+1},0,\dots$ nor $0, a_j,0,a_{j+1}\dots$.
Then, $(4.2)$ holds true by taking again the proof of Proposition 4.1. 
When $k\ge l+1$, we have
$$
t_k= \cases
{-m + q_{k-1} & if $k= l+2j, \quad (j\ge 1)$,\cr
-m+q_k & if $k=l+2j+1,\quad (j\ge 0),$\cr}
$$
and 
$$
\sigma_k= \cases
{-\theta_{k-1} & if $k= l+2j, \quad (j\ge 1)$,\cr
-\theta_k & if $k=l+2j+1,\quad (j\ge 0).$\cr}
$$

Assume first that $k$ has the same parity as $l$, namely $k= l+2j$ for some $j\ge 1$. Then $\sigma_k = -\theta_{k-1}$ and $t_k = -m+q_{k-1}$. In order to check $(4.2)$, we distinguish three subcases.
Assume first $ q\le q_{k}-1$. Then $\| q\theta \| \ge | \theta_{k-1} | $ with equality only when $q=q_{k-1}$. If $q\not= q_{k-1}$, we have 
$$
\| q\theta \| > |\theta_{k-1} | ,
$$
so that $q\theta$ and $q\theta -\theta_{k-1}$ are located in the same open interval of length one, so that $(4.2)$ holds true. If $q=q_{k-1}$, then we have
$$
q_{k-1}\theta -\theta_{k-1} = p_{k-1} \quad \hbox{and} \quad q_{k-1}\theta = p_{k-1}+ \theta_{k-1},
$$
so that $(4.2)$ holds, since $\theta_{k-1}$ is positive, noting that $k-1 = l-1+2j$ is even. 
Assume secondly that $q=q_k$. Then, 
$$
q_{k}\theta -\theta_{k-1} = p_{k} +\theta_k-\theta_{k-1} \quad \hbox{and} \quad q_{k}\theta = p_{k}+ \theta_{k}. 
$$
This shows that $(4.2)$ holds, since both numbers $ \theta_k$ and $\theta_k -\theta_{k-1}$ 
are negative with absolute value less than $1$. 
Assume thirdly that $q= q_k +v$ for some integer $v$ with $1 \le v \le t_k = -m + q_{k-1}$. Then,
$$
q\theta -\theta_{k-1} = p_{k} +\theta_k-\theta_{k-1}+v\theta \quad \hbox{and} 
\quad q\theta = p_{k}+ \theta_{k} +v\theta.
$$
Notice now that  $\| v\theta \| \ge | \theta_{k-2} | $ with equality 
only when $v=q_{k-2}$. If $v\not= q_{k-2}$, we have 
$$
\| v\theta \| > |\theta_{k-2} | . 
$$
Then, $q\theta$ and $q\theta -\theta_{k-1}$ are located in the same open interval of length one, since
$$
| \theta_{k-2} |  \ge | \theta_k  |  + | \theta_{k-1}| , 
$$
so that $(4.2)$ holds true. If $v=q_{k-2}$, then we have
$$
q\theta -\theta_{k-1} = p_k+ p_{k-2} + \theta_k -\theta_{k-1} +\theta_{k-2} \quad \hbox{and} \quad q\theta =p_k +  p_{k-2}+ \theta_{k}+ \theta_{k-2},
$$
so that $(4.2)$ holds, since $\theta_k -\theta_{k-1} +\theta_{k-2}$ and $ \theta_{k}+ \theta_{k-2}$ are both negative with absolute value less than  $1$. 

We assume now  that $k= l+2j+1 $ for some $j\ge 0$. Then $\sigma_k = -\theta_{k}$ and $t_k = -m+q_{k}$. We distinguish again three subcases.
Assume first that $q\le q_k-1$. Then, $\| q\theta \| \ge | \theta_{k-1} |$, so that $q\theta$ and $q\theta -\theta_{k}$ are located in the same open interval of length one. It folllows that $(4.2)$ holds. 
Assume secondly that $q=q_k$. Then, 
$$
q_{k}\theta -\theta_{k} = p_{k}  \quad \hbox{and} \quad q_{k}\theta = p_{k}+ \theta_{k}.
$$
This shows that $(4.2)$ holds, since  $ \theta_k$ is positive because $k = l+2j+1$ is even. 
Assume thirdly that $q= q_k +v$ for some integer $v$ with $1 \le v \le t_k = -m + q_{k}$. Then,
$$
q\theta -\theta_{k} = p_{k}  +v\theta \quad \hbox{and} \quad q\theta = p_{k}+ \theta_{k} +v\theta. 
$$
Notice now that  $\| v\theta \| \ge | \theta_{k-1} |  > | \theta_k | $. 
Thus $(4.2)$ holds. All cases have been checked. 

When $l$ is even, the numbers $\theta_{l+ 2j}$ (resp. $\theta_{l+2j+1}$)  
turn to be positive (resp. negative), and the above argumentation remains  
valid provided that we replace the usual integer part $\lf \cdot \rf$ by the upper integer part $\lc \cdot \rc$.

\cqfd

To illustrate this statement, take $\rho = 0$, $a_1 = 5$, $a_2 = 3$, $a_3 = 2$; then 
$$
V_0 = 0, \quad V_1 = 1 0^4, \quad V_2 = 1 0^4 1 0^4 1 0^4 0, \quad 
$$
$$
V_3 = 1 0^4 1 0^4 1 0^4 1 0^4 0 1 0^4 1 0^4 1 0^4 0 = 1 0^4 1 0^4 1 0^4 1 0^5 1 0^4 1 0^4 1 0^5,
$$
$$
V'_0 = 0, \quad V'_1 =  0^4 1, \quad V'_2 = 0 0^4 1 0^4 1 0^4 1, \quad 
$$
$$
V'_3 = 0 0^4 1 0^4 1 0^4 1 0 0^4 1 0^4 1 0^4 1 0^4 1
= 0^5 1 0^4 1 0^4 1 0^5 1 0^4 1 0^4 1 0^4 1.
$$
Note also that 
$$
M_0 = 0, \quad M_1 = 0^4 1, \quad M_2 = 0^4 1 0^4 1 0^4 1 0, 
$$
$$
M_3 = 0^4 1 0^4 1 0^4 1 0 0^4 1 0^4 1 0^4 1 0 0^4 1 = 0^4 1 0^4 1 0^4 1 0^5 1 0^4 1 0^4 1 0^5 1. 
$$
By induction, we check that $V_n$ is the mirror image of $V'_n$. 
We know that $M_n^{--}$ is a palindrome. We also have that 
$$
{}^{-}V_n^{-} = {}^{-}(V'_n)^{-} = M_n^{--},
$$
where ${}^{-}W$ means the word $W$ deprived of its first letter.  
In other words, for $n \ge 1$, the 
words $V_n$ and $V'_n$ deprived of their first and last letters are equal to the palindrome $M_n^{--}$.

\vskip 5mm

\centerline{\bf 5. Repetitions in a Sturmian word}

\vskip5mm

We keep our notation. 
Recall that $\bx$ denotes an arbitrary Sturmian word of slope $\theta$.

We show that Proposition 1 of \cite{BKLN} can be deduced 
from the recursion formulae for the words $V_k$ 
and we give further  informations on the occurrence of the various cases. 
Proposition 5.1 will be used in the next section to compare $\bx$ with 
four families of (shifted for two of them) periodic words,  
depending on a parameter $k$, constructing  thus families 
of strong rational approximations to the associated Sturmian number.   

\proclaim Proposition 5.1.     
Let $k$ be an integer with $k \ge 0$.  
Then, there exist a uniquely determined    
non-empty suffix $U_k$ of $M_k M_{k+1} = (M_k)^{a_{k+1}+1} M_{k-1} $
and an integer $\ta_{k+1}$  
 such that 
$$
\ta_{k+1} \in \{a_{k+1}, a_{k+1} + 1 \}  
$$
and 
$$
{\bx} = U_k (M_k)^{\ta_{k+1}} M_{k-1} M_{k}^{-} \dots  
$$ 
More precisely, when $a_{k+2}-b_{k+2}\ge 2$, we have 
$$
U_k = R_{k+1} \quad \hbox{and} \quad \ta_{k+1} = a_{k+1}.
$$
When $a_{k+2}-b_{k+2}=  1$, we have 
$$
U_k = R_{k+1} \quad \hbox{and} \quad \ta_{k+1} =
\cases{ a_{k+1}+1 &  if \  $b_{k+3} < a_{k+3}$, \cr
a_{k+1} & if  \ $b_{k+3}= a_{k+3}$.
\cr}
$$
When $a_{k+2}= b_{k+2}$, we have $U_k = R_{k}M_{k+1} $. Moreover 
$\ta_{k+1} = a_{k+1}$, unless
$$ 
a_{k+2}= 1, \, a_{k+3}-b_{k+3}\ge 2,
$$
or
$$
a_{k+2}=1, \, a_{k+3}= 1, \, b_{k+3}= 0, \, a_{k+4}= b_{k+4}, 
$$
in which cases $\ta_{k+1}= a_{k+1}+ 1$.

\medskip

 {\bf Remark.} The fact that ${\bx} = U_k (M_k)^{\ta_{k+1}} M_{k-1} M_{k}^{-} \dots $ means that 
after the prefix of length $|U_k|$, we have exactly $\ta_{k+1} + 1$ copies of $M_k$, followed 
by the prefix of $M_k$ of length $q_{k-1} - 2$,
since $M_{k-1}M_k^-$ and $M_kM_{k-1}^-$
differ only by their last letter. 
In addition, we observe that when $a_{k+2}= b_{k+2}$ we have $b_{k+1} = 0$ 
and we take $U_{k-1} = R_k$.  

\medskip

\pro
 The idea of the proof is to show that  the prefix of $\bx$ of length $2q_{k+1} + q_k-1$ coincides with  one of the three words $V_{k+1}^2V_k^{-}$ or   $V_{k+1}V_kV_{k+1}^{-}$ or  $V_kV_{k+1}V_{k+1}^{-}$.
 
 
 

Assume first that $a_{k+2} -b_{k+2}\ge 2$. Then
$$
V_{k+2} = V_{k+1}^{a_{k+2}-b_{k+2}} V_k V_{k+1}^{b_{k+2}}= V_{k+1}^2V_k^{-}\dots , 
$$
observing that $V_k^{-} $ is a prefix of $V_{k+1}$ 
(this follows from Definition 2.2).  
But 
$$
\eqalign{
V_{k+1}V_{k+1}V_k&= R_{k+1}T_{k+1}R_{k+1}T_{k+1}R_kT_k \cr
& = R_{k+1}M_{k+1}T_{k+1}R_kT_k  \cr
& = R_{k+1}M_{k+1}M_k^{b_{k+1}}T_kR_kT_k  \cr
&= R_{k+1}M_{k+1}M_k^{b_{k+1}+1}T_k = R_{k+1} M_k^{a_{k+1}}M_{k-1}M_k\dots   
\cr}
$$

Assume secondly that $a_{k+2} -b_{k+2}=1$ and that $a_{k+3}-b_{k+3}\ge 1$. Then
$$
V_{k+3} = V_{k+2}^{a_{k+3}-b_{k+3}}V_{k+1}V_{k+2}^{b_{k+3}}=V_{k+2} V_{k+1}^{-}\dots= (V_{k+1} V_k V_{k+1}^{b_{k+2}+1})^{-}\dots = V_{k+1}V_{k}V_{k+1}^{-}\dots
$$
Actually, we can be more precise and claim that $V_{k+1}V_{k}V_{k+1} V_k^{-}$ is a prefix of ${\bx}$. 
This is obvious unless $b_{k+2} = 0$  (then $a_{k+2}=1$ and $V_{k+2} = V_{k+1}V_k$) and $a_{k+3} - b_{k+3} = 1$ and $b_{k+3} = 0$ (then $V_{k+3} = V_{k+2}V_{k+1}= V_{k+1}V_{k}V_{k+1}$). 
Assume that these three equalities hold. If $a_{k+4}> b_{k+4}$, we have
$$
V_{k+4} = V_{k+3} V_{k+2}^{-} \ldots = V_{k+1} V_{k}V_{k+1} V_{k+2}^{-} \ldots = V_{k+1} V_k  V_{k+1}V_{k+1} V_k^{-} \ldots , 
$$
then $V_{k+1} V_k  V_{k+1} V_k^{-}$ is indeed a prefix of $\bx$. Otherwise, we have
$$
V_{k+4} = V_{k+2} V_{k+3} \ldots = V_{k+1}V_k V_{k+1} V_{k} V_{k+1}\ldots = V_{k+1} V_k  V_{k+1} V_k \ldots ,
$$
and the same conclusion holds. 

We claim that 
$$
V_{k+1}V_{k}V_{k+1}V_k= R_{k+1}M_k^{a_{k+1}+1}M_{k-1}M_k^{b_{k+1}+1}T_k,
$$
which yields that  $U_k = R_{k+1} $ and $\ta_{k+1}= a_{k+1}+1$. 
For the proof,  we distinguish two cases, either $a_{k+1}> b_{k+1}$, or $a_{k+1}= b_{k+1}$. In the first case, we have
$$
T_{k+1}= M_k^{b_{k+1}}T_k \quad {\rm and} \quad R_{k+1}= R_kM_k^{a_{k+1}-b_{k+1}-1} M_{k-1},
$$
so that we compute
$$
\eqalign{
V_{k+1}V_{k}V_{k+1} V_k&= R_{k+1}T_{k+1}R_{k}T_{k}R_{k+1}T_{k+1} R_k T_k\cr
&= R_{k+1}M_k^{b_{k+1}}T_k R_k T_k R_kM_k^{a_{k+1}-b_{k+1}-1}M_{k-1}M_k^{b_{k+1}}T_{k} R_kT_k \cr  
& =R_{k+1}M_k^{a_{k+1}+1}M_{k-1}M_k^{b_{k+1}+1}T_{k}.
\cr}
$$
For the latter case, we have 
$$
T_k= T_{k-1}, \,\, R_k = R_{k-1}M_{k-1}^{a_k-1}M_{k-2}, \,\, T_{k+1}= M_k^{a_{k+1}}T_{k-1}= M_k^{a_{k+1}}T_k, \,\,  R_{k+1}= R_{k-1}.
$$
Thus
$$
\eqalign{
V_{k+1}V_{k}V_{k+1} V_k&= R_{k+1}T_{k+1}R_{k}T_{k}R_{k+1}T_{k+1}R_kT_k
\cr
&= R_{k+1}M_k^{a_{k+1}}T_{k-1}R_{k-1}M_{k-1}^{a_k-1}M_{k-2} T_{k-1} R_{k-1} M_k^{a_{k+1}} T_kR_kT_k  
\cr  
&= R_{k+1}M_k^{a_{k+1}+1} M_{k-1} M_k^{a_{k+1}+1}T_k.
\cr}
$$
The claim is established.

Assume thirdly that $a_{k+2} -b_{k+2}=1$ and that $a_{k+3}=b_{k+3}$. 
Then,  $b_{k+2}= 0$ and $a_{k+2}=1$. We find
$$
V_{k+3}  = V_{k+1}V_{k+2}^{a_{k+3} }  
= V_{k+1}(V_{k+1}V_k)^{a_{k+3}} = V_{k+1}^2 V_k \dots.
$$
The first case shows that $U_k = R_{k+1}$ and $\ta_{k+1}= a_{k+1}$,  as asserted.

Suppose finally that $a_{k+2}=b_{k+2}$. Then $b_{k+1}= 0$ 
and $a_{k+3}> b_{k+3}$, since $a_{k+3}= b_{k+3}$ should yield $a_{k+2}= b_{k+2}=0$. Thus,
$$
V_{k+3}= V_{k+2}^{a_{k+3}-b_{k+3}}V_{k+1}V_{k+2}^{b_{k+3}} 
=V_{k+2}V_{k+1}^{-}\dots = V_kV_{k+1}^{a_{k+2}}V_{k+1}^{-}\dots = V_kV_{k+1}V_{k+1}^{-}\dots
$$
Here, again, we can be more precise and show that $\bx$ is either of the form 
$$
 \bx = V_kV_{k+1}V_kV_{k+1}V_k^{-}\dots , \eqno{(5.1)}
 $$ 
 or  of the form
 $$
 \bx = V_k V_{k+1}V_{k+1} V_k^{-}\dots.\eqno{(5.2)}
 $$ 
 If $a_{k+2}\ge 2$, we have
 $$
 V_{k+3}= V_{k+2} V_{k+1}^{-} \dots = V_kV_{k+1}^{a_{k+2}}V_k^{-}\dots = V_kV_{k+1}^2 V_k^{-}\dots
 $$
 Thus $(5.2)$ holds. When $a_{k+2}=1$ and $a_{k+3}-b_{k+3}\ge 2$, we have $V_{k+2}= V_kV_{k+1}$ and 
 $$
V_{k+3}= V_{k+2}^2V_{k+1}^{-}\dots= V_kV_{k+1}V_kV_{k+1}V_k^{-}\dots
$$ 
Thus $(5.1)$ holds. When $a_{k+2} = 1$, $a_{k+3} -b_{k+3} = 1$ and $b_{k+3}\ge 1$, we have
$$
V_{k+3}= V_{k+2}V_{k+1}V_{k+2}^{b_{k+3}}= V_{k}V_{k+1}^2V_{k+2}^{b_{k+3}},
$$
so that $(5.2)$ holds true. When $a_{k+2}=1$, $a_{k+3}= 1$ and $b_{k+3}=0$, we have $V_{k+3} = V_kV_{k+1}^2$. If $a_{k+4}-b_{k+4}\ge 1$, we have 
$$
V_{k+4}= V_{k+3} V_{k+2}^{-} \dots = V_kV_{k+1}^2 V_k^{-}\dots
$$
so that $(5.2)$ holds, while 
$$
V_{k+4} = V_{k+2}V_{k+3}\dots = V_kV_{k+1}V_kV_{k+1}^2\dots
$$
if $a_{k+4}= b_{k+4}$. Then $(5.1)$ holds. Now, we compute
$$
\eqalign{
V_{k}V_{k+1}V_{k+1} V_k&= R_{k}T_{k}R_{k+1}T_{k+1}R_{k+1}T_{k+1}R_kT_k= R_{k}T_k(R_k M_k^{a_{k+1}-1}M_{k-1}) M_{k+1}  T_{k} R_kT_k 
\cr  
& 
= R_{k}M_k^{a_{k+1}}M_{k-1}M_k^{a_{k+1}}M_{k-1}M_kT_k
=R_k M_{k+1} M_k^{a_{k+1}}M_{k-1} T_k,
\cr}
$$
and
$$ 
\eqalign{
V_{k}V_{k+1}V_{k} V_{k+1}V_k &= R_{k}T_{k}R_{k+1}T_{k+1}R_{k}T_{k}R_{k+1}T_{k+1} R_{k}T_{k}
\cr  & = R_kT_kR_k M_k^{a_{k+1}-1}M_{k-1}T_kR_kT_kR_kM_k^{a_{k+1}-1}M_{k-1}T_kR_kT_k
\cr
& = R_k M_k^{a_{k+1}}M_{k-1}M_k^{a_{k+1}+1}M_{k-1}M_kT_k
\cr
& = R_kM_{k+1}M_k^{a_{k+1}+1}M_{k-1}M_kT_k.
\cr}
$$
Thus $U_k= R_kM_{k+1}$ in both cases. We have  $\ta_{k+1}= a_{k+1}$ when $(5.2)$ holds, while $\ta_{k+1}= a_{k+1}+1$ whenever $(5.1)$ is satisfied.
\cqfd

\bigskip


We have used at several places   
the obvious property that $V_k^{-}$ is a prefix of $V_{k+1}$, which holds since by definition 
$V_k$, $V_{k+1}$ and $\bx$ share the same prefix of length $q_k-1$. 
A question which arises naturally is to know when  $V_k$ is a prefix of $V_{k+1}$.

\proclaim Proposition 5.2. For any $k\ge 0$, the word $V_k$ is a prefix of $V_{k+1}$
if and only if the sequence $b_1, \dots, b_{k+1}$ differs from $0,a_2, 0, a_4, \dots, a_{k+1}$ when $k$ is odd, 
or differs from $ a_1-1,0,a_3, 0,\dots , a_{k+1}$ when $k$ is even. 

For $k \ge 0$, let $W_k$ denote the longest common prefix of $V_{k+1} V_k$ and $V_k V_{k+1}$. 

\proclaim Lemma 5.3. 
We have $W_0 = 0^{a_1-1-b_1}$ and $W_{k+1} = V_{k+1}^{a_{k+2} - b_{k+2}} W_k$ for $k \ge 0$.    
Consequently,  the length $w_k$ of $W_k$ is given by 
$$
w_k= a_1-1-b_1 + \sum_{j=1}^k(a_{j+1}-b_{j+1})q_j = q_{k+1}+q_k -t_{k+1}-2, \quad k \ge 0. 
$$

\pro 
Recall that $V_{-1} = 1$, $V_0 = 0$, and $V_1 = V_0^{a_1-1-b_1} V_{-1} V_0^{b_1}$. 
This implies that $V_0 V_1 = V_0^{a_1-b_1} V_{-1} V_0^{b_1}$, thus
$$
W_0 = V_0^{a_1-1-b_1}, \quad w_0 = a_1-1-b_1.
$$
We proceed by induction. Let $k \ge 0$ be an integer. 

Assume first that $a_{k+2} - b_{k+2} \ge 1$. 

Since $V_{k+2} = V_{k+1}^{a_{k+2} - b_{k+2}} V_k V_{k+1}^{b_{k+2}}$, we get 
$$
V_{k+2} V_{k+1} = V_{k+1}^{a_{k+2} - b_{k+2}} V_k V_{k+1}^{b_{k+2}+1}
$$
and 
$$
V_{k+1} V_{k+2} = V_{k+1}^{a_{k+2} - b_{k+2} + 1} V_k V_{k+1}^{b_{k+2}},
$$
thus
$$
W_{k+1} = V_{k+1}^{a_{k+2} - b_{k+2}} W_k.
$$

Assume now that $a_{k+2} = b_{k+2}$. In that case, we know that $b_{k+1} = 0$. 
Then,
$$
V_{k+2} V_{k+1} = V_k V_{k+1}^{a_{k+2}  + 1} = V_k V_k^{a_{k+1}} V_{k-1} V_k^{a_{k+1}} V_{k-1}
V_{k+1}^{a_{k+2} - 1}
$$
and
$$
V_{k+1} V_{k+2} = V_k^{a_{k+1}} V_{k-1} V_k V_{k+1}^{a_{k+2}},
$$
thus
$$
W_{k+1} = V_k^{a_{k+1}} W_{k-1} = W_k = V_{k+1}^{a_{k+2} - b_{k+2}} W_k.
$$
Since $q_k$ is the length of $V_k$, this proves the lemma. 
\cqfd

\noindent {\it Proof of Proposition 5.2. }      
The word $V_k$ is a prefix of $V_{k+1}$ exactly when $w_k\ge q_k$.  
Lemma 5.3 tells us that $w_k \ge q_k$ if and only if $t_{k+1}\le q_{k+1}-2$. 
Observe finally that $t_{k+1}\le q_{k+1}-1$ with equality if and only if 
$$
b_{k+1}= a_{k+1},\,  b_k = 0, \, b_{k-1}= a_{k-1}, \, b_{k-2}= 0, \dots.
$$
This completes the proof. 
\cqfd

\vskip 5mm

\centerline{\bf 6. The sequence of convergents contributing to the exponent of irrationality}

\vskip 5mm

In this section and the next one, $b \ge 2$ is an integer and $\xi$ denotes one of the numbers 
$\xi_b (\theta, \rho)$ or $\xi'_b (\theta, \rho)$. 
We analyze the convergents which contribute to the exponent of irrationality of $\xi$, which we 
call `strong convergents'. 
According to \cite{BuKim19}, all of them are obtained by truncating the $b$-ary expansion of $\xi$ 
and completing by periodicity. Thus, their denominators are either of the form $b^s - 1$ (purely 
periodic case) or $b^r (b^s - 1)$ (existence of a preperiod). 

We adopt the following conventions of writing. Any finite word $Y=y_1 \ldots y_r$
with letters in $\{ 0, \dots, b-1\}$ is as well viewed as the natural integer
$$
Y = y_1b^{r-1}+ \cdots + y_r,
$$
whose sequence of $b$-ary digits is given by $Y$.     
Then, for any words $Y=y_1\ldots y_r$ and $Z=z_1 \ldots z_s$, we have the $b$-ary expansions 
$$
{ Z \over b^s -1} =0. Z^\infty,
$$
and
$$
 { Y Z - Y \over b^r(b^s-1)}= 0. Y Z^\infty,
$$ 
where $0.z_1z_2\dots ={z_1\over b} + {z_2\over b^2}+ \cdots$  and  $YZ$ stands 
for the number whose $b$-ary  sequence of digits is the concatenation 
of the words $Y$ and $Z$, that is, $YZ = y_1 \ldots y_r \, z_1\ldots z_s$.

Let  $u$ and $v$ be two positive quantities depending upon a parameter $k$. 
As usual, we write $u\asymp v$  when there exist positive constants $c_1$ and $c_2$,
independent of $k$,   
such that $c_1 u \le v \le c_2 u$.

The candidates for the sequence of strong convergents belong to four types.  
We label them by the index $k\ge 0$.    
The sequences of finite words $(R_k)_{k \ge 0}, (T_k)_{k \ge 0}, (V_k)_{k \ge 0}$ are   
given by Theorem 2.1 (or Theorem 4.2), but we  now  replace 
the alphabet $\{0, 1\}$ by $\{0, b-1\}$.      
We recall that $R_0= 0$ and $T_0$ is the empty word.   
Below, the height  means the logarithmic height $\log H/\log b$, that is, roughly speaking, the largest   
exponent of $b$ appearing in the denominator.

The first possible convergent is
$$
(1)_k =   { R_{k+1}-R_k\over b^{r_k}( b^{r_{k+1}-r_k}-1)},
$$
with height $\asymp r_{k+1}$ and $b$-ary expansion
$$
(1)_k= 0.R_k(M_k^{a_{k+1}-b_{k+1}-1} M_{k-1})^{\infty}. 
$$  
Of course $(1)_k$ is meaningful only when $r_{k+1} > r_k$, that is to say when $a_{k+1}- b_{k+1}\ge 1$. 
The second candidate is
$$
(2)_k =   { R_{k+1}T_k\over  b^{r_{k+1}+ t_k}-1},
$$
with height $\asymp r_{k+1}+ t_k$, associated to the periodic word 
$(R_{k+1}T_k)^\infty$. The third is
$$
(3)_k =   { R_{k+1}M_k-M_k\over  b^{r_{k+1}}(b^{q_k}-1)},
$$
with height $\asymp r_{k+1}+ q_k$, associated to the word $R_{k+1}M_k^\infty$.  The fourth is 
$$
(4)_k =   { V_{k+1}\over  b^{q_{k+1}}-1},
$$
with height $\asymp q_{k+1}= r_{k+1}+t_{k+1}$,
associated to the periodic word $V_{k+1}^\infty=(R_{k+1}T_{k+1})^\infty$. 

We say that a rational $x$ precedes another one $y$, and we write $x \prec y$, 
when the height of $x$ is less than the height of $y$.
Clearly 
$$
(1)_k \prec (2)_k \preceq (4)_k.
$$
We have $(2)_k = (4)_k$ exactly when $t_k=t_{k+1}$, that is to say when $b_{k+1}=0$. Then,
$$
(1)_k \prec (2)_k = (4)_k \prec  (3)_k.
$$
If $b_{k+1}\ge 1$, we have 
$$
t_{k+1}\ge t_k+q_k,
$$
so that
$$
(1)_k \prec (2)_k \prec (3)_k \prec (4)_k,
$$
in this case.  When $a_{k+1}=b_{k+1}$, obviously $b_{k+1} \ge 1$, so that the above inequality 
$$
(2)_k \prec (3)_k \prec (4)_k
$$
hold with $(1)_k$ being omitted.

An important observation is that we have the following coincidences between levels $k-2$, $k-1$ and $k$.

If $a_{k+1}-b_{k+1}= 1$, then we have 
$$
(1)_k= (3)_{k-1},
$$
since 
$$
R_k(M_k^{a_{k+1}-b_{k+1}-1} M_{k-1})^{\infty}= R_kM_{k-1}^\infty.
$$

If $a_{k+1}= b_{k+1}$, then we have
$$
(2)_k = (4)_{k-2}, 
$$
since we have $R_{k+1}= R_{k-1}$ and $T_k= T_{k-1}$ (because $b_k=0$), so that 
$$
R_{k+1}T_k = R_{k-1}T_{k-1} = V_{k-1}.
$$
If, in addition, $b_{k-1} = 0$, then $T_{k-1} = T_{k-2}$ and $(2)_{k-2} = (2)_k = (4)_{k-2}$.  

Observe also that if $a_{k+1}-b_{k+1}=1$, then we have
$$
(2)_k= {V_kV_{k-1}\over b^{q_k+q_{k-1}}-1},
$$
since 
$$
R_{k+1}T_k = R_kM_{k-1}M_{k-1}^{b_k}T_{k-1}=R_kT_kR_{k-1}T_{k-1}= V_kV_{k-1}, 
$$
noting that
$$
T_kR_{k-1} =T_{k-1}V_{k-1}^{b_k}R_{k-1}= T_{k-1}(R_{k-1}T_{k-1})^{b_k}R_{k-1} 
= (T_{k-1}R_{k-1})^{b_k+1}=M_{k-1}^{b_k+1}.
$$

To go further for linking consecutive blocks (with indices  $k-1$ and $k$),  
we need to know when the rationals $(1)_k, (2)_k, (3)_k, (4)_k$ are indeed convergents. 
We indicate as well in the next proposition the value of the exponential rate  
of approximation $\mu_k(j)$ such 
$$
| \xi -(j)_k | \asymp {1 \over H((j)_k)^{\mu_k(j)}} = {1 \over b^{h((j)_k)\mu_k(j)}},
$$
for all large $k$ and $1 \le j \le 4$, where $h((j)_k)$ is the base-$b$ logarithm of
the height  $H((j)_k)$ of $(j)_k$. 
We determine in 
which cases the exponent $\mu_k(j)$ is bigger than $2$, thanks to Proposition 5.1. 
To that purpose, let us introduce the following quantities
$$
\nu_k(1) = 1+ {r_{k+1}+t_k\over r_{k+1}}, \quad \nu_k(2)
= 1+ {r_{k+1}+ q_k\over r_{k+1}+ t_k}, 
$$
$$
\nu_k(3) = 1+ { q_{k+1}\over r_{k+1}+ q_k},
\quad \nu_k(4)= 1+ { r_{k+2}\over q_{k+1}}.
$$
They are equal  to one plus the ratio of   the height of two consecutive points
in the sequence $ \dots (1)_k, (2)_k,(3)_k, (4)_k, (1)_{k+1},\dots$.
Then, we can state the following criterion.

\proclaim
Proposition 6.1. 
Let $k\ge 2$ be an integer such that $t_{k-1}$ is positive.   
The rational
$$
(1)_k = { R_{k+1}-R_k\over b^{r_k} (b^{r_{k+1} - r_k}-1)}
$$
is a convergent to $\xi$  if and only if
 $$
 a_{k+1}-b_{k+1}\ge 1,  \, a_{k+2}-b_{k+2} \ge 1 \quad  \hbox{and then} \quad  \mu_k(1)= \nu_k(1),   
 $$
 or
 $$
b_k \ge 1, \,   a_{k+1}= 1, \, b_{k+1}= 0,  \, a_{k+2}=b_{k+2}  \quad \hbox{and then} \quad \mu_k(1)= \nu_{k-1}(3).
$$
 The rational
$$
(2)_k = { R_{k+1}T_k\over b^{r_{k+1}+t_k}- 1}
$$
is a convergent to $\xi$ if and only if
$a_{k+2}-b_{k+2} \ge 1$ and then
$$
  \mu_k(2) = \cases{  \nu_k(2) & if $ b_{k+1}\ge 1$, \cr
  \nu_k(4) & if $ b_{k+1}= 0, \, a_{k+3}-b_{k+3}\ge 1$,  \cr
   \nu_{k+2}(2) & if $b_{k+1}=0, \, a_{k+3}=b_{k+3}$  .
  \cr}    
 $$
  The rational
$$
(3)_k = { R_{k+1}M_k-R_{k+1}\over b^{r_{k+1}}(b^{q_k}- 1)}
$$
is a convergent to $\xi$ if and only if
$$
b_{k+1}\ge1,\,   a_{k+2}-b_{k+2} \ge 2 \quad  \hbox{and then} \quad  \mu_k(3)= \nu_{k}(3),        
$$
or 
$$
  a_{k+2}-b_{k+2} = 1,\,  a_{k+3}-b_{k+3} \ge 1  
  \quad  \hbox{and then} \quad  \mu_k(3)= \nu_{k+1}(1),        
$$
or
$$
b_{k+1}\ge1,\,    a_{k+2}=1, \,  b_{k+2} =0,\, 
a_{k+3}= b_{k+3}      \quad  \hbox{and then} \quad  \mu_k(3)= \nu_{k}(3),     
$$
 The rational
$$
(4)_k = { V_{k+1}\over b^{q_{k+1}}- 1}
$$
is a convergent to $\xi$ if and only if
$$
 a_{k+2}-b_{k+2} \ge 2,\,  a_{k+3}-b_{k+3} \ge 1  
  \quad  \hbox{and then} \quad  \mu_k(4)= \nu_{k}(4),     
 $$
 or
 $$
b_{k+1}=0  ,  \,  a_{k+2}-b_{k+2} =  1, \,  a_{k+3}-b_{k+3} \ge  1  
\quad  \hbox{and then} \quad  \mu_k(4)= \nu_{k}(2)= \nu_k(4),     
$$
or  
$$
 a_{k+3}= b_{k+3} 
 \quad  \hbox{and then} \quad  \mu_k(4)= \nu_{k+2}(2)= 1+ \nu_k(4).    
 $$

\pro
We  only prove Proposition 6.1 assuming that  $t_{k-1}$ is large enough.  In fact, crude estimates 
of the constants involved in the symbols $\asymp$ show that the lower bound $b^{t_{k-1}} \ge 4$ 
is sufficient for our purpose. Relaxing the assumption to  $t_{k-1}\ge 1$ 
follows from an alternative argumentation  which will be given in the next Section 7. Our present approach is based on Legendre's theorem 
asserting that $P/Q$ is a convergent to $\xi$ when  $| \xi -P/Q | < 1/ (2Q^2)$.      

Let $\bx$ be the Sturmian word composed of the $b$-ary digits of $\xi$. 

For $(1)_k$, 
the relevant assumption is $a_{k+1} -b_{k+1} \ge 1$. Assume first that $ a_{k+2}-b_{k+2} \ge 1$. Then, Proposition 5.1  gives
$$    
\eqalign{
\bx   &=  R_{k+1} M_{k}^{\ta_{k+1}} M_{k-1} M_{k}^{-} \ldots   \cr
&= R_{k} M_{k}^{a_{k+1} - b_{k+1} - 1} M_{k-1} M_{k}^{\ta_{k+1}} M_{k -1}\ldots  \cr
& = R_{k} M_{k}^{a_{k+1} - b_{k+1} - 1} M_{k-1} M_{k}^{a_{k+1} - b_{k+1} - 1} 
M_{k}^{\ta_{k+1} -a_{k+1} + b_{k+1}+1} M_{k-1} \ldots  \cr
& = R_{k} M_{k}^{a_{k+1} - b_{k+1} - 1} M_{k-1} M_{k}^{a_{k+1} - b_{k+1}- 1} 
M_{k} M_{k-1} \ldots ,
}
$$
to be compared with the word $ R_{k} (M_{k}^{a_{k+1} - b_{k+1} - 1} M_{k-1} )^{\infty}$. When $a_{k+1}-b_{k+1}\ge 2$, we can write
$$
 R_{k} (M_{k}^{a_{k+1} - b_{k+1} - 1} M_{k-1} )^{\infty}= 
 R_{k} M_{k}^{a_{k+1} - b_{k+1} - 1} M_{k-1} M_{k}^{a_{k+1} - b_{k+1}- 1} 
M_{k-1} M_{k}\dots 
 $$
 to  obtain the estimate
$$
|  \xi - (1)_k | 
 \asymp {1 \over b^{r_{k} + 2(a_{k+1}- b_{k+1} - 1)q_{k} + 2 q_{k -1}+ q_{k}}}  
 = {1 \over b^{2 r_{k+1} + t_{k}}}.  
$$
When $a_{k+1}-b_{k+1} = 1$ the same estimate holds, since  then
$$
\x = R_kM_{k-1}M_kM_{k-1} \dots = R_kM_{k-1}^{a_k+1}M_{k-2}M_{k-1}\dots, 
$$
while 
$$
 R_{k} (M_{k}^{a_{k+1} - b_{k+1} - 1} M_{k-1} )^{\infty}=R_kM_{k-1}^\infty =R_kM_{k-1}^{a_k+1} M_{k-1} M_{k-2}\dots.
 $$
Thus $(1)_k$ is a convergent to $\xi$ and  
$$
\mu_k(1)= {2 r_{k+1} +t_k \over r_{k+1} }= 1+ {r_{k+1}+ t_k \over r_{k+1}}= \nu_k(1).  
$$
When $a_{k+2}= b_{k+2}$, we have $b_{k+1}=0, U_k = R_kM_{k+1}$, and Proposition 5.1 gives
$$
\bx   =  U_k \ldots   
= R_{k} M_{k}^{a_{k+1}  - 1} M_k M_{k-1} \ldots  
$$
We distinguish two subcases. If $a_{k+1}\ge 2$, we write
$$
R_k(M_k^{a_{k+1}-1}M_{k-1})^\infty = R_kM_k^{a_{k+1}-1}M_{k-1}M_k M_k^{a_{k+1}-2} M_{k-1}\dots.
$$

Thus,
$$
|  \xi  -(1)_k |   
 \asymp {1 \over b^{r_{k} + a_{k+1} q_{k} +  q_{k -1}}}  
 = {1 \over b^{ r_{k+1} + q_{k} }},  
$$
so that
$$
(r_{k+1}+q_k)-2r_{k+1}= q_k-r_{k+1}= q_k-(r_k+ (a_{k+1}-1)q_k +q_{k-1}) = t_k-q_{k+1} +q_k
$$
is negative, since $q_{k+1}> 2 q_k$. Therefore $(1)_k$ is not a convergent in this subcase. 
When $a_{k+1}=1$, write
$$
\bx =R_kM_kM_{k-1}\dots   = R_kM_{k-1}^{a_k}M_{k-2}M_{k-1}\dots, 
$$
while 
$$
R_kM_{k-1}^\infty = R_k M_{k-1}^{a_k} M_{k-1}M_{k-2} \dots
$$
Thus,
$$
| \xi - (1)_k | 
 \asymp {1 \over b^{r_{k} + (a_k+1)q_{k-1} +  q_{k -2}}}  
 = {1 \over b^{ r_{k} + q_{k}+q_{k-1} }}
  = {1 \over b^{ r_{k+1} +q_{k} }},  
$$
so that
$$
\eqalign{
(r_{k+1}+q_{k})-2r_{k+1}= -r_{k+1}+q_{k} & =-(r_k+q_{k-1})+q_k \cr
& = t_k-q_{k-1}
=t_{k-1}+ (b_k-1) q_{k-1}. \cr
}   
$$
We conclude by noticing that $t_{k-1}+(b_k-1)q_{k-1}$ is positive if $b_k\ge 1$ and negative when $b_k=0$.
Thus, 
$$
\mu_k(1)= {r_{k+1} +q_k\over r_{k+1}} = 1 + {q_k\over r_{k+1} } = 1 + {q_k\over r_k +q_{k-1}} = \nu_{k-1}(3).
$$
Observe that, in this case, we have the ordering
$$
(1)_k = (3)_{k-1} \prec  (2)_{k+1}=(4)_{k-1},
$$
while $(2)_k, (3)_k, (4)_k$ and $(1)_{k+1}$ are not convergents to $\xi$.

We now deal with  $(2)_k$.  Assume first that $a_{k+2} -b_{k+2}\ge 1$.  

In the subcase   $a_{k+1}-b_{k+1}\ge 1$ and $b_{k+1}\ge 1$, 
we have $R_{k+1}= R_kM_k^{a_{k+1}-b_{k+1}-1}M_{k-1}$ and Proposition 5.1 gives
$$
\eqalign{
\bx   &=  R_{k+1} M_{k}^{\ta_{k+1}} M_{k-1} M_{k}^{-}\ldots   \cr
&= R_{k} M_{k}^{a_{k+1} - b_{k+1} - 1} M_{k-1} M_{k}^{\ta_{k+1}} M_{k -1}M_k^{-}\ldots  \cr
& = R_{k} M_{k}^{a_{k+1} - b_{k+1} - 1} M_{k-1} M_{k}^{a_{k+1} - b_{k+1} } 
M_{k}^{\ta_{k+1} -a_{k+1} + b_{k+1}} M_{k-1}M_k^{-} \ldots  \cr
& = R_{k} M_{k}^{a_{k+1} - b_{k+1} - 1} M_{k-1} M_{k}^{a_{k+1} - b_{k+1}} 
M_{k} M_{k-1} \ldots 
 }
$$
since $\ta_{k+1}-a_{k+1}+b_{k+1}\ge 1$. 
Comparing with the word 
$$ 
(R_{k} M_{k}^{a_{k+1} - b_{k+1} - 1} M_{k-1} T_{k})^{\infty} 
= R_k M_{k}^{a_{k+1} - b_{k+1} - 1} M_{k-1} M_{k}^{a_{k+1} - b_{k+1} } M_{k-1}M_k\dots , 
$$
we obtain
$$
  | \xi - (2)_k | 
  \asymp {1\over b^{r_k +2(a_{k+1}-b_{k+1})q_k  +2q_{k-1}}}
 = {1 \over b^{2(r_{k+1} + t_{k}) + r_{k}}}.
$$
Thus, $(2)_k$ is a convergent of $\xi$  and
$$
\mu_k(2) ={  2(r_{k+1} +t_k) + r_k \over r_{k+1}+t_k }= 2 + {r_k\over r_{k+1}+t_k}  = 1+{r_{k+1}+q_k\over r_{k+1}+t_k}= \nu_k(2).
$$ 

In the subcase  $a_{k+1}=b_{k+1}$ (and thus $b_{k+1} \ge 1$), we have  $(2)_k = (4)_{k-2}$.  
Assuming temporarily that Proposition 6.1 has been checked for $(4)_{k-2}$, it yields that $(2)_k$ is again  a convergent to $\xi$ with exponent $\mu_k(2)=\mu_{k-2}(4) = \nu_k(2)$ as asserted.

Consider finally the subcase  $b_{k+1}=0$. 
Then $(2)_k= (4)_k$ and  Proposition 6.1 for $(4)_k$, tells us  that $(2)_k$ is indeed a convergent to $\xi$ with exponent $\mu_k(2)= \mu_k(4)$ which will be computed below.

Assume now that $a_{k+2}= b_{k+2}$. Then, Proposition 5.1  gives
$$
\bx = R_k M_{k+1}\dots= R_kM_k^{a_{k+1}}M_{k-1}\dots =  
R_kM_k^{a_{k+1}-1}M_kM_{k-1}\dots , 
$$
while
$$
(R_{k+1}T_k)^\infty = R_kM_k^{a_{k+1}-1}M_{k-1}T_kR_k\dots
=R_kM_k^{a_{k+1}-1}M_{k-1}M_k\dots
$$
since $b_{k+1}=0$. It follows that
$$
\vert \xi -(2)_k \vert \asymp {1 \over b^{ r_k + a_{k+1}q_k +q_{k-1}}}
= {1 \over b^{ r_{k+1}+ q_k}} ={1\over b^{ 2(r_{k+1}+t_k)- (q_{k+1}-r_k)}} .
$$  
Then $(2)_k $ is not a convergent to $\xi$.

We now deal with $(3)_k$. Assume first that $a_{k+2}-b_{k+2}\ge 1$. 
 Proposition 5.1 gives
$$
\bx   =  R_{k+1} M_{k}^{\ta_{k+1}} M_{k-1} M_{k}^{-} \ldots   
  $$ 
Since
$$
R_{k+1}M_k^\infty =  R_{k+1} M_{k}^{\ta_{k+1}} M_kM_{k -1}\dots   
$$
we obtain the estimate
$$
|  \xi- (3)_k |   
\asymp {1\over b^{r_{k+1}+ (\ta_{k+1}+1)q_k +q_{k-1}}}.     
$$ 
Write
$$
\eqalign{
(r_{k+1}+ (\ta_{k+1}+1)q_k +q_{k-1}) -2 (r_{k+1}+q_{k})& =-r_{k+1}+ (\ta_{k+1}-1)q_k +q_{k-1} 
\cr
&= t_{k+1}  + (\ta_{k+1}-a_{k+1}-1 )q_{k} .
\cr
}
$$
If $a_{k+2}-b_{k+2}=1$ and $a_{k+3}-b_{k+3}\ge 1$, we know that $\ta_{k+1}= a_{k+1}+1$, 
so that $ t_{k+1}  + (\ta_{k+1}-a_{k+1}-1 )q_{k}>0$. 
If $a_{k+2}-b_{k+2}\ge 2$, or if $a_{k+2}=1, b_{k+2}=0, a_{k+3}=b_{k+3}$, 
we know that $\ta_{k+1}= a_{k+1}$, so that 
$$
 t_{k+1}  + (\ta_{k+1}-a_{k+1}-1 )q_{k}= t_k + (b_{k+1}-1)q_k.
 $$
 Now, $t_k + (b_{k+1}-1)q_k$ is positive when $b_{k+1}\ge 1$ and negative when $b_{k+1}=0$.  We get the three cases announced.  
 Concerning  the exponent $\mu_k(3)$, we find 
 $$
 \mu_k(3) = { r_{k+1}+ (\ta_{k+1}+1)q_k+ q_{k-1}\over r_{k+1}+q_k}.     
 $$
 When $\ta_{k+1}=a_{k+1}$, we get 
 $$
 \mu_k(3) = { r_{k+1}+ q_{k+1}+q_k\over r_{k+1}+ q_k} =1+ {q_{k+1}\over r_{k+1}+ q_k} = \nu_k(3),
 $$
 while, in the case  $\ta_{k+1}=a_{k+1} + 1$, we have  
 $$
 \mu_k(3) = { r_{k+1}+ q_{k+1}+2q_k\over r_{k+1}+ q_k} =1+ {q_{k+1}+q_k\over r_{k+1}+ q_k} 
 =1+{r_{k+2} +t_{k+1}\over r_{k+2}}= \nu_{k+1}(1),
 $$
 since $r_{k+2} = r_{k+1}+ q_k$ when $a_{k+2}-b_{k+2}=1$.
 It remains for us to prove that $(3)_k$ is not a convergent when $a_{k+2}= b_{k+2}$. 
 Then, $b_{k+1}= 0$ and $r_{k+1}= r_{k} + (a_{k+1}-1)q_k +q_{k-1}$. In this case, Proposition 5.1 gives
 $$
 \bx = R_kM_{k+1}\dots = R_k M_k^{a_{k+1}-1}M_k M_{k-1}\dots 
 $$
 while
 $$
 R_{k+1}M_k^\infty = R_kM_k^{a_{k+1}-1} M_{k-1} M_k\dots.  
 $$
 Thus 
 $$
 \vert \xi -(3)_k \vert \asymp{ 1\over b^{ r_k+ a_{k+1}q_k +q_{k-1}} }
 \asymp {1\over b^{r_{k+1}+ q_k}}, 
 $$
and  $(3)_k$ is not a convergent to $\xi$. 
 
 For the last  rational
 $$
 (4)_k = { V_{k+1}\over b^{q_{k+1}}-1}
 $$
Proposition 5.1 tells us that  $\bx = V_{k+1}^2V_k^-\dots $ whenever 
$$
 a_{k+2}-b_{k+2} \ge 2  
 $$
or
$$
a_{k+2}= 1 \quad \hbox{and} \quad b_{k+2} =  0
\quad \hbox{and} \quad a_{k+3}= b_{k+3}  .
$$
Then, the initial exponent of repetition of $V_{k+1}$ is clearly larger than $2$, 
so that $(4)_k$ is a convergent to $\xi$.
When $a_{k+2}-b_{k+2}=1$ and $a_{k+3}-b_{k+3} \ge 1$, we have 
$$
\bx = V_{k+1}V_kV_{k+1}^-\dots 
$$
By Lemma 5.3, the common prefix $W_k$ to $V_kV_{k+1}$ and $V_{k+1}V_k$ has length    
$$
w_k= a_1-1-b_1 + \sum_{j=1}^k(a_{j+1}-b_{j+1})q_j = q_{k+1}+q_k -t_{k+1}-2.
$$
Noting that
$$
t_{k+1} = \sum_{j=0}^k b_{j+1}q_j
$$
is larger or smaller than $q_k$ when $b_{k+1}\ge 1$ or $b_{k+1}=0$, we deduce that $(4)_k$ 
is then a convergent to $\xi$ when $b_{k+1}= 0$ and is not when $b_{k+1}\ge 1$. This yields the case 
 $$
b_{k+1}=0  \quad \hbox{and} \quad a_{k+2}-b_{k+2} =  1 \quad \hbox{and} \quad a_{k+3}-b_{k+3} \ge  1. 
$$
When $ a_{k+3}-b_{k+3}\ge 1$ and $a_{k+2}-b_{k+2}\ge 1$, Proposition 5.1, 
with $k$ replaced by $k+1$, tells us that 
$$  
\bx  = R_{k+2}M_{k+1}\dots = R_{k+1}M_{k+1}^{a_{k+2}-b_{k+2}-1}M_kM_{k+1} \dots ,
$$
while 
$$
V_{k+1}^\infty = R_{k+1}M_{k+1}^\infty = R_{k+1}M_{k+1}^{a_{k+2}-b_{k+2}-1}M_{k+1}M_{k}\dots.  
$$
It follows that
$$
| \xi -(4)_k  | \asymp {1 \over b^{ r_{k+1} +(a_{k+2}-b_{k+2})q_{k+1}+q_k} }= {1\over b^{r_{k+2} +q_{k+1}}}.
$$
Thus, 
$$
\mu_k(4) = {r_{k+2} + q_{k+1}\over q_{k+1}} = 1 + {r_{k+2} \over q_{k+1}} = \nu_k(4).
$$
Notice that $ \nu_k(4) = \nu_k(2)$ in the case $ a_{k+2}-b_{k+2}= 1$ and $b_{k+1} = 0$, 
since $(1)_{k+1}= (3)_k$ and $(4)_k=(2)_k$. 

When $ a_{k+3}=b_{k+3}$, Proposition 5.1 with $k$ replaced by $k+1$,  gives 
$$  
\bx  =R_{k+1}M_{k+2}M_{k+1}^{\ta_{k+2}}\dots= R_{k+1}
M_{k+1}^{a_{k+2}}M_kM_{k+1}\dots 
$$
It follows that
$$
| \xi -(4)_k  | \asymp {1 \over b^{ r_{k+1} +(a_{k+2}+1)q_{k+1}+q_k} }= {1\over b^{r_{k+2} +2q_{k+1}}},
$$
since $r_{k+2} = r_{k+1}+(a_{k+2}-1)q_{k+1}+ q_k$. Thus, 
$$
\mu_k(4) = {r_{k+2} +2 q_{k+1}\over q_{k+1}} = 2 + {r_{k+2} \over q_{k+1}} =1+  \nu_k(4) = \nu_{k+2}(2),
$$
noting that  
$$
r_{k+2} +q_{k+1}=  r_{k+1}+q_{k+2}= r_{k+3}+ q_{k+2},
$$
and
$$
q_{k+1}= r_{k+1}+ t_{k+1}= r_{k+3}+ t_{k+2},
$$
since $b_{k+2}=0$.

When $a_{k+2}=b_{k+2}$, the word $\bx$ has a prefix of  the form
$$
\bx= V_{k+2}\dots = V_kV_{k+1}^{b_{k+2}}\dots
$$
and the common prefix  of $V_{k+1}^\infty$ and $\bx$  has length at most 
$$
w_k \le q_{k+1}+ q_k -2 < 2 q_{k+1}.
$$
Thus, $(4)_k$ cannot be  a convergent to $\xi$. 
\cqfd

\bigskip

The next proposition describes a tail of the sequence of strong convergents ordered by increasing height. 
We start with the cyclic sequence ${\cal S}$
$$
(1)_0, (2)_0, (3)_0, (4)_0, (1)_1, \dots , (4)_{k-1}, (1)_k, (2)_k, (3)_k, (4)_k, (1)_{k+1}, \dots
$$
built with the $(j)_k$. As already observed,  some elements of ${\cal S}$ 
may coincide and the height function is not necessarily increasing along ${\cal S}$.  
Assume that $\bx$ differs from $\bc_\theta$, so that $t_k$ is positive for any $k\ge h$ and some $h\ge 1$. 
Then, let  ${\cal S}^+$ be the tail of ${\cal S}$ formed by the elements $(j)_k$ with  $k\ge h+1$. 
Assuming moreover  that $1 \le b_k \le a_k-2$ for every $k\ge h+1$,   
Proposition 6.1 tells us that the sequence of strong convergents $(j)_k$, 
restricted to the indices $k\ge h+1$,  coincides with ${\cal S}^+$.
Otherwise, the following modifications are needed.

\bigskip

\noindent{ \bf Proposition 6.2.} 
{\it A tail of the ordered sequence of strong convergents to $\xi$ is obtained  
by applying to ${\cal S}^+$ the following replacement rules.

(i) Assume $a_{k+2}=b_{k+2}$. When $b_k\ge 1$,   we replace the string of seven elements 
$(4)_{k-1}, \dots, (2)_{k+1}$ by the single element $ (4)_{k-1}= (2)_{k+1}$.
  When $b_k=0$,  we replace the string of  nine elements $(2)_{k-1}, \dots , (2)_{k+1}$ 
  by the single element $(2)_{k-1}=(4)_{k-1}= (2)_{k+1}$. 

(ii) Assume $a_{k+2}-b_{k+2}= 1$ and $a_{k+3}-b_{k+3}\ge 1$. 
When $b_{k+1} \ge 1$, 
we replace the three elements $(3)_k, (4)_k,(1)_{k+1}$ by the single element
$ (3)_k= (1)_{k+1}$,  and the four elements $(2)_k, (3)_k, (4)_k,(1)_{k+1}$ by the pair 
$$ 
(2)_k =(4)_k\prec (3)_{k}= (1)_{k+1},
$$ 
when $b_{k+1}= 0$.

(iii) Assume that $a_{k+2}-b_{k+2}\ge 2$  and $a_{k+3}-b_{k+3}\ge 1$. When $b_{k+1}=0$, we replace the three elements $(2)_k, (3)_k, (4)_k$ by  the single element $ (2)_k =(4)_k$.
}

\bigskip
\noindent
{\bf Remark.}  Observe that there is no overlap for the above replacement rules, since the case $(i)$ cannot appear for two consecutive indices $k$ by Ostrowski's numeration rules.

\bigskip

\pro We check  in each case $(i)$, $(ii)$ and $(iii)$ that the elements $(j)_k$ in  ${\cal S}$   
which are erased do not belong to the list provided by Proposition 6.1,  
while the remaining ones belong indeed to the list.

For instance, in the case $(i)$ with $b_k=0$, Proposition 6.1 tells us that 
$\mu_{k-1}(2)= \nu_{k+1}(2)$. Moreover,
$$
 (2)_{k-1}=(4)_{k-1}= (2)_{k+1} \prec (3)_{k+1}
 $$
 are convergents to $\xi$, while the intermediate rationals 
 $ (3)_{k-1}, \, (1)_k, \, (2)_k, \, (3)_k, \, (4)_k, \, (1)_{k+1}$ are not, 
 as can be verified by reading  the necessary and sufficient conditions 
 displayed in Proposition 6.1 for each element involved.  
\cqfd

\bigskip

It will be proved in Proposition 7.2  
that the subset  of   convergents to $\xi$ given by Proposition 6.1  
provides all the convergents contributing to the   irrationality  exponent of $\xi$. 
We thus obtain the 
\proclaim 
Theorem 6.3. 
The irrationality exponent of $\xi$ is equal to
$$
\limsup_{k \to + \infty} \, 
\max\{ \mu_k(1), \mu_k(2), \mu_k(3), \mu_k(4)\}
=\max\{ \nu(1), \nu(2), \nu(3),\nu(4)\},
$$
where 
$$
\eqalign{
\nu(1) &= \limsup_{k \to + \infty} \{ \nu_k(1) \, : \, a_{k+1}-b_{k+1} \ge 1 
\quad{\rm and} \quad a_{k+2}-b_{k+2} \ge 1\},
\cr
\nu(2) &= \limsup_{k \to + \infty} \{ \nu_k(2) \, : \, a_{k+2}-b_{k+2} \ge 1 \},
\cr
\nu(3) &= \limsup_{k \to + \infty} \{ \nu_k(3)   \},
\cr
\nu(4) &= \limsup_{k \to + \infty} \{ \nu_k(4)  \}.
\cr}
$$

\pro  
For any convergent $(j)_k$ to $\xi$, we have expressed  $\mu_k(j)$ 
as some value $ \nu_{k'}(j')$, thanks to Proposition 6.1. Conversely,
for any given $\nu_k(j)$,  we analyze   
under which conditions it contributes 
to the exponent of irrationality of $\xi$. For instance, Proposition 6.1 tells us that $\nu_k(1)$ 
occurs exactly when $a_{k+1}-b_{k+1} \ge 1$ and $a_{k+2}-b_{k+2}\ge 1$, 
leading to the definition of $\nu(1)$. Similarly, $\nu_k(2)$ appears in Proposition 6.1  if and only if 
$$ 
b_{k+1} \ge 1, \ \ a_{k+2}-b_{k+2}\ge 1,
$$
or
$$
b_{k+1}=0, \ \ a_{k+2}-b_{k+2} = 1, \ \ a_{k+3}-b_{k+3} \ge 1,
$$
or 
$$
 a_{k+1}=b_{k+1}.
$$
Remark first that the third case   
is included in the first case, because $a_{k+1} = b_{k+1}$ implies $a_{k+2}-b_{k+2}\ge 1$ 
by Ostrowski's rules. Recall  that $(2)_k=(4)_k$ when $b_{k+1} =0 $. 
Observe now that the assumptions $b_{k+1} =0$ and $a_{k+2}-b_{k+2}\ge 1$ 
yield the inequality $\nu_k(2) \le \nu_k(4)$, with equality if and only if $a_{k+2} -b_{k+2} =1$, 
since $t_{k+1}=t_k$ and
$$
\nu_k(2)= 1+{r_{k+1} +q_k \over r_{k+1}+t_k}=1+{ r_{k+2}- (a_{k+2}-b_{k+2}-1)q_{k+1}\over q_{k+1} }
\le 1+{r_{k+2}\over q_{k+1}} = \nu_k(4).
$$
 We may thus remove the condition $b_{k+1}\ge 1$ in the first case, since the additional contributions  
 are taken into account  by $\nu(4)$. Finally, the single constraint $a_{k+2}-b_{k+2} \ge 1$ remains. 
 We are thus led to introduce the quantity $\nu(2)$. 
 
 We now deal with the contribution of $\nu_k(4)$. It  occurs in Proposition 6.1 exactly when
 $$
 a_{k+2}-b_{k+2} \ge 2 , \ \ a_{k+3} -b_{ k+3} \ge 1,
 $$
 or 
 $$
 b_{k+1}= 0, \ \ a_{k+2}-b_{k+2}= 1, \ \ a_{k+3}-b_{k+3}\ge 1.  
 $$
 Observe that  $\nu_k(4) =1+{r_{k+2}\over q_{k+1}} $ 
 is at most equal to $2$ when $b_{k+1}\ge 1$ and $a_{k+2}-b_{k+2}=1$, since then
 $$
 r_{k+2} = r_{k+1}+ q_k = \cases{r_k+q_{k+1} -b_{k+1}q_k \le q_{k+1}-t_k & if $a_{k+1}-b_{k+1}\ge 1$, \cr
 r_{k-1} +q_k \le q_{k+1} & if $ a_{k+1}=b_{k+1}$.\cr}
 $$
 We may thus forget the condition $b_{k+1}=0$ in the second case above. 
 Observe also that $\nu_k(4) < 2$ when $a_{k+2}=b_{k+2}$. 
 It remains the constraint $ a_{k+3}-b_{k+3}\ge 1$. Notice however  
 that we may remove this last constraint as asserted. Indeed, when $a_{k+3}=b_{k+3}$, 
 Proposition 6.1 tells us that $(4)_k = (2)_{k+2}$ is a convergent to $\xi$ 
 with approximation exponent $\nu_{k+2}(2) = 1 + \nu_k(4)$. 
 Since $a_{k+4}-b_{k+4}\ge 1$ by Ostrowski's rules, the number $\nu_{k+2}(2)> \nu_k(4)$ 
 is taken into account by $\nu(2)$. We may thus  define $\nu(4)$ unconditionally as above.
 
We finally deal with the contribution of $\nu_k(3)$. It appears when
$$
b_{k+1} \ge 1, \ \  a_{k+2}-b_{k+2} \ge 2,
 $$
 or 
 $$
 b_{k+1}\ge  1, \ \ a_{k+2} = 1, \ \ b_{k+2}=0, \ \ a_{k+3}=b_{k+3}.
 $$
We may relax the constraints as follows. We first forget the assumption $b_{k+1} \ge 1$, since when $b_{k+1} =0$, we have $r_{k+1}= r_k +q_{k+1}-q_k$, so that 
$$
\nu_k(3) = 1 + {q_{k+1} \over r_{k+1}+ q_k} = 1 + {q_{k+1} \over r_k + q_{k+1}} < 2. 
$$
We may also relax  the assumptions 
$a_{k+2}=1, b_{k+2}=0,a_{k+3}=b_{k+3}$ in the second case above to $a_{k+2}-b_{k+2}=1$, since  
when $a_{k+2}-b_{k+2}= 1$ and $a_{k+3}-b_{k+3} \ge 1$, we have $(3)_k = (1)_{k+1}$, while
$$
\nu_k(3) = 1 + {q_{k+1} \over r_{k+1} +q_k} = 1 + {q_{k+1}\over r_{k+2}} 
< 1+ {q_{k+1} + q_k\over r_{k+2}} = 1 + {r_{k+2} +t_{k+1}\over r_{k+2}}= \nu_{k+1}(1).
$$
The additional contributions are then covered by $\nu(1)$. 
It remains the constraint $a_{k+2}-b_{k+2} \ge 1$. But when $a_{k+2}= b_{k+2}$, we have $b_{k+1}=0$, so that $\nu_k(3) \le 2$, as already observed. 
 \cqfd



\vskip 5mm

\goodbreak

\centerline{\bf 7. The partial quotients}

\vskip 5mm

We keep the notation of the previous section. 

For $k \ge 0$, recall that we have set 
$$
c_k = b^{r_{k} + q_{k-1}} \, 
{b^{(a_{k+1} - b_{k+1} - 1 ) q_{k}} - 1 \over b^{q_{k}} - 1},
\quad
d_k = b^{t_{k}} - 1,
$$
$$
e_k = b^{r_{k}} - 1, \quad 
f_k = b^{t_{k}} \, {b^{b_{k+1} q_{k}} - 1 \over b^{q_{k}} - 1}.
$$
The integers $c_k, d_k, e_k, f_k$ are positive, unless $t_k = 0$ (and then $d_k = 0$)   
or $b_{k+1} = 0$ (and then $f_k = 0$) or $a_{k+1} - b_{k+1} \le 1$ (and then $c_k = 0$ if 
$a_{k+1} - b_{k+1} = 1$, while otherwise $c_k = - b^{r_{k-1}} = - e_{k-1} - 1$, by (2.3)).

Recall that we have defined the possible convergents by 
$$
(1)_k =   { R_{k+1}-R_k\over b^{r_k}( b^{r_{k+1}-r_k}-1)}, 
\quad
(2)_k =   { R_{k+1}T_k\over  b^{r_{k+1}+ t_k}-1}, \quad k \ge 0, 
\quad 
$$
$$
(3)_k =   { R_{k+1}M_k-M_k\over  b^{r_{k+1}}(b^{q_k}-1)},
\quad
(4)_k =   { V_{k+1}\over  b^{q_{k+1}}-1}, \quad k \ge 0.
$$
Put also 
$$
(4)_{-1} =   { V_{0}\over  b^{q_{0}}-1} =   { {0}\over  b -1}. 
$$
From a Diophantine point of view, 
$(1)_k$ is meaningful only when $r_{k+1} > r_k$, that is to say when $a_{k+1}- b_{k+1}\ge 1$. 
Nevertheless, it can be formally defined as well when $r_{k+1} < r_k$, in which case 
numerator and denominator are negative integers.

In the sequel, the notation $(1)_{k} = c_k \cdot (4)_{k-1} + (3)_{k-1}$ means that 
the numerator (resp., denominator) of $(1)_{k}$ is equal to 
$c_k$ times the numerator (resp., denominator) of $(4)_{k-1}$ plus the 
numerator (resp., denominator) of $(3)_{k-1}$. 
With some abuse of notation,  
$$
(2)_k-(1)_k = { R_{k+1}T_{k}-(R_{k+1}-R_k) \over b^{r_{k+1}+t_k} -1-(b^{r_{k+1}}-b^{r_k})}
$$ 
stands below 
for the ratio of the difference between the numerators and denominators of $(2)_k$ and $(1)_k$.

\proclaim Lemma 7.1. 
For $k \ge 0$, we have the following relations: 
$$
(1)_{k} = c_k \cdot (4)_{k-1} + (3)_{k-1},  \quad (k \not= 0), 
$$
$$
(2)_{k} - (1)_k = d_k   \cdot (1)_{k} + (4)_{k-1},  
$$
$$
(2)_{k} = 1 \cdot ( (2)_{k} - (1)_{k} ) +  (1)_{k},
$$
$$
(3)_{k} = e_k  \cdot (2)_{k} + ( (2)_{k} - (1)_{k} ), 
$$
$$
(4)_{k} = f_k \cdot (3)_{k} + (2)_{k}. 
$$

\pro
Let us begin with the first equality. If $a_{k+1} - b_{k+1} \ge 1$, then
$$
\eqalign{
c_k (b^{q_k} - 1) + b^{r_k} (b^{q_{k-1}} - 1) 
& = b^{r_{k} + q_{k-1}} \, 
(b^{(a_{k+1} - b_{k+1} - 1 ) q_k} - 1 ) + b^{r_k} (b^{q_{k-1}} - 1) \cr
& = b^{r_{k}  +( a_{k+1} - b_{k+1} - 1) q_k  + q_{k-1} } - b^{r_k} \cr
& = b^{r_{k+1}} - b^{r_k}, \cr
}
$$
which is the denominator of $(1)_k$. 
Likewise, we have
$$
\eqalign{ 
& V_k \times b^{r_{k} + q_{k-1}} \, 
{b^{(a_{k+1} - b_{k+1} - 1 ) q_{k}} - 1 \over b^{q_{k}} - 1} \cr
&= R_k T_k \times (b^{r_{k} + q_{k-1}} + b^{r_{k} + q_{k-1} + q_{k}} 
+ \ldots + b^{r_{k} + q_{k-1} + (a_{k+1} - b_{k+1} - 2 ) q_{k}} ) \cr
& = (R_{k} T_{k})^{a_{k+1} - b_{k+1} - 1} b^{r_{k} + q_{k-1}} \cr
& = (R_{k} T_{k})^{a_{k+1} - b_{k+1}  - 1} R_{k} M_{k-1} - R_{k} M_{k-1}  \cr
& = R_{k+1}  - R_{k} M_{k-1} = ( R_{k+1}  - R_{k} ) - (R_{k} M_{k-1}  - R_k), \cr
}
$$
if $a_{k+1} - b_{k+1} \ge 2$, while 
$$
V_k \times b^{r_{k} + q_{k-1}} \, 
{b^{(a_{k+1} - b_{k+1} - 1 ) q_{k}} - 1 \over b^{q_{k}} - 1} = 0 = ( R_{k+1}  - R_{k} ) - (R_{k} M_{k-1}  - R_k),
$$
if $a_{k+1} - b_{k+1} = 1$, because we then have $R_{k+1} = R_k M_{k-1}$. 
In both cases we end up with the numerator of $(1)_k$ minus the numerator of $(4)_{k-1}$.   

Now, assume that $a_{k+1} = b_{k+1}$. Then, $c_k = - b^{r_{k-1}}$ and we check that 
$$
( - b^{r_{k-1}}) (b^{q_k} - 1) + b^{r_k} (b^{q_{k-1}} - 1) = b^{r_{k-1}} - b^{r_k} = b^{r_{k+1}} - b^{r_k},
$$
since $r_{k-1} + q_k = r_k + q_{k-1}$ and $r_{k-1} = r_{k+1}$. 
As for the numerators, we have 
$$
\eqalign{
V_k \times ( - b^{r_{k-1}}) & = - V_k R_{k-1} + R_{k-1}  \cr
& = - R_k T_k R_{k-1} + R_{k-1} \cr
& = - R_k M_{k-1} + R_{k+1} = ( R_{k+1}  - R_{k} ) - (R_{k} M_{k-1}  - R_k), \cr
}
$$
which confirms our claim. 

For the second equality, observe that 
$$
b^{t_k} \bigl( b^{r_k}( b^{r_{k+1}-r_k}-1) \bigr) + (b^{q_k} - 1) 
= b^{q_k + r_{k+1}} - b^{t_k + r_k} + b^{q_k} - 1 = b^{q_k + r_{k+1}} -  1
$$
is the denominator of $(2)_k$. Note also that 
$$
b^{t_k} (R_{k+1} - R_k) = R_{k+1} T_{k} - R_k T_{k} = R_{k+1} T_{k} - V_k
$$
is the numerator of $(2)_k$ minus the numerator of $(4)_k$. This completes 
the proof of the second equality. The third one is a tautology. 
The remaining two equalities are proved in a similar way than the second one. 
We omit the details. 
\cqfd

Define two sequences $(P_j)_{j \ge -1}$ and $(Q_j)_{j \ge -1}$ of integers by setting 
$$
P_{-1} = b-1, \ \ Q_{-1}= 0, \ \  
P_0 = 0, \ \ Q_0 = b-1 , 
$$
and, denoting by 
$(\alpha_j)_{j \ge 1}$ the sequence of integers $c_0, d_0, 1, e_0, f_0, c_1, \ldots$, 
$$
P_{j+2} = \alpha_{j+2} P_{j+1} + P_j, \quad 
Q_{j+2} = \alpha_{j+2} Q_{j+1} + Q_j, \quad j \ge -1.
$$
Since 
$$
c_0= {b^{a_1-b_1} -b\over b-1}, \ \ d_0= 0, \ \   e_0 = b-1,
$$
we get
$$
\displaylines{
P_1 = b-1,\ \ Q_1= b^{a_1-b_1}-b, \ \ P_2 = P_0, \ \ Q_2 =Q_0, 
\cr 
\ \ P_3 = b-1, \ \ Q_3 = b^{a_1-b_1}-1, \ \ P_4 = (b-1)^2, \ \ Q_4 = b^{a_1-b_1}(b-1), \dots .
\cr}
$$ 
Thus 
$$
{P_1\over Q_1}=(1)_0, \ \ {P_2\over Q_2}= (2)_0 -  (1)_0,\ \  
{P_3\over Q_3}= (2)_0, \ \ {P_4\over Q_4}= (3)_0,  \ \ \dots.
$$
Using Lemma 7.1, we check by induction on $j \ge 1$ that the greatest 
prime divisor of the integers $P_j$ and $Q_j$ is equal to $b-1$ and 
that $P_j$ and $Q_j$ are the numerator and 
denominator of a fraction of one of the five types
$(1)_{k}$, $(2)_{k} - (1)_k$, $(2)_k$, $(3)_k$, $(4)_k$, more precisely, they correspond to 
\smallskip
$*$ the fraction $(1)_k$ if $\alpha_j = c_k$; 

$*$ the fraction $(2)_k - (1)_k$ if $\alpha_j = d_k$; 

$*$ the fraction $(2)_k$ if $\alpha_j = 1$; 

$*$ the fraction $(3)_k$ if $\alpha_j = e_k$; 

$*$ the fraction $(4)_k$ if $\alpha_j = f_k$.
\smallskip



We explain below how to derive the sequence 
of partial quotients of $\xi$ from the sequence $(\alpha_j)_{j \ge 1}$. 

To do this, we work with matrices and recall that      
$$
\eqalign{
\pmatrix{
P_{j+1} & P_j \cr
Q_{j+1} & Q_j \cr
}
& =
\pmatrix{
P_{j} & P_{j-1} \cr
Q_{j} & Q_{j-1} \cr
}
\cdot
\pmatrix{
\alpha_{j+1} & 1 \cr
1 & 0 \cr
} \cr
& =\pmatrix{
0 & b-1 \cr
b-1 & 0 \cr
}    
\pmatrix{
\alpha_{1} & 1 \cr
1 & 0 \cr
}
\cdots 
\pmatrix{
\alpha_{j+1} & 1 \cr
1 & 0 \cr
}, 
\quad j \ge 0. 
\cr
}
$$
So we have a product of elementary integer $2$ by $2$ matrices $\pmatrix{
\alpha_j  & 1 \cr
1 & 0 \cr
}$, exactly as in the continued 
fraction algorithm. Here, however, some coefficients $\alpha_j$ may be $0$ or negative. 
The point is that it is possible to transform this formal infinite   
product into a product of elementary integer $2$ by $2$ matrices $\pmatrix{
\alpha''_j  & 1 \cr
1 & 0 \cr
}$ where all the $\alpha''_j$'s are positive. This defines a regular continued fraction and 
we show that this is precisely the continued fraction expansion of~$\xi$.

Simple calculations show that for nonnegative integers $x$ and $y$ we have
$$
\pmatrix{
x  & 1 \cr
1 & 0 \cr
}
\cdot
\pmatrix{
0 & 1 \cr
1 & 0 \cr
}
\cdot 
\pmatrix{
-x-1 & 1 \cr
1 & 0 \cr
}
=
\pmatrix{
-1 & 1 \cr
1 & 0 \cr
} 
\eqno (7.1)
$$
and
$$
\pmatrix{
y & 1 \cr
1 & 0 \cr
}
\cdot 
\pmatrix{
1 & 1 \cr
1 & 0 \cr
}
\cdot
\pmatrix{
-1  & 1 \cr
1 & 0 \cr
}
\cdot
\pmatrix{
y & 1 \cr
1 & 0 \cr
}
\cdot 
\pmatrix{
1 & 1 \cr
1 & 0 \cr
}
=
\pmatrix{
0 & 1 \cr
1 & 1 \cr
}.
\eqno (7.2)
$$
If for some integer $j\ge 1$ the integer $\alpha_{j+5} = c_{k+1}$  
is negative, then $c_{k+1} = - e_k - 1$ and,     
as $b_{k+1}= 0$, we get $d_{k+1} = d_k$, $f_k = 0$,      
and the septuple $(\alpha_{j+1}, \ldots , \alpha_{j+7})$
is equal to $(d_k, 1, e_k, 0, -e_k-1, d_k, 1)$. Consequently, by (7.1) and (7.2), we have 
$$
\pmatrix{
\alpha_{j+1} & 1 \cr
1 & 0 \cr
} 
\cdots
\pmatrix{
\alpha_{j+7} & 1 \cr
1 & 0 \cr
} 
= 
\pmatrix{
0 & 1 \cr
1 & 1 \cr
}.
$$
We derive that
$$
\eqalign{
\pmatrix{
P_{j+8} & P_{j+7} \cr
Q_{j+8} & Q_{j +7} \cr
}
& = 
\pmatrix{
P_{j-1} & P_{j-2} \cr
Q_{j-1} & Q_{j-2} \cr
}
\cdot
\pmatrix{
\alpha_{j} & 1 \cr
1 & 0 \cr
} 
\cdots
\pmatrix{
\alpha_{j+8} & 1 \cr
1 & 0 \cr
} \cr
& = 
\pmatrix{
P_{j-1} & P_{j-2} \cr
Q_{j-1} & Q_{j-2} \cr
}
\cdot 
\pmatrix{
\alpha_{j} & 1 \cr
1 & 0 \cr
} 
\cdot 
\pmatrix{
0 & 1 \cr
1 & 1 \cr
}
\cdot
\pmatrix{
\alpha_{j+8} & 1 \cr
1 & 0 \cr
} 
\cr
& = 
\pmatrix{
P_{j-1} & P_{j-2} \cr
Q_{j-1} & Q_{j-2} \cr
}
\cdot
\pmatrix{
\alpha_j + \alpha_{j+8} + 1 & 1 \cr
1 & 0 \cr
} \cr
& = 
\pmatrix{
(\alpha_j + \alpha_{j+8} + 1) P_{j-1} + P_{j-2} & P_{j-2} \cr
(\alpha_j + \alpha_{j+8} + 1) Q_{j-1} + Q_{j-2}  & Q_{j -2} \cr
} \cr
& = 
\pmatrix{
(c_k + e_{k+1} + 1) P_{j-1} + P_{j-2} & P_{j-1} \cr
(c_k + e_{k+1} + 1) Q_{j-1} + Q_{j-2}  & Q_{j -1} \cr      
}. \cr
}   \eqno (7.3)
$$
This shows that $P_{j-1}$ is followed by $P_{j+8} = (c_k + e_{k+1} + 1) P_{j-1} + P_{j-2}$, and   
similarly for $Q_{j-1}$. 

Consider now the sequence $(\alpha'_j)_{j \ge 1}$ constructed inductively from 
$(\alpha_j)_{j \ge 1}$ as follows. We put $\alpha'_j = \alpha_j$ for $j < j_0$, where 
$j_0\ge 1$ is the smallest integer such that $\alpha_{j_0} = c_k$, with $\alpha_{j_0 +5}=c_{k+1} < 0$.  
Then, we put $\alpha'_{j_0} = c_k + e_{k+1} + 1$ and $\alpha'_{j_0 + 1} = \alpha_{j_0 + 9} = f_{k+1}$. 
We continue with $\alpha'_{j_0 + 2} = c_{k+2}$,
 unless $c_{k+3} < 0$, in which case 
we put $\alpha'_{j_0 + 2} = c_{k+2} + e_{k+3} + 1$. 
And so on. 
The sequence $(\alpha'_j)_{j \ge 1}$ is well-defined since $c_k$ and $c_{k+1}$ cannot be 
simultaneously negative. 

Said differently, for each index $k$ such that $c_{k+1} < 0$, we replace the $10$ 
consecutive partial quotients $c_k, d_k, \ldots , e_{k+1}, f_{k+1}$ by the $2$ partial 
quotients $c_k + e_{k+1} + 1, f_{k+1}$. Let us add that $f_{k+1}$ is positive since $b_{k+2}$ is 
positive. 

We have constructed from $(\alpha_j)_{j \ge 1}$ a sequence of nonnegative integers $(\alpha'_j)_{j \ge 1}$. 
Define 
$$
P'_{-1} = b-1, \ \ P'_0= 0, \ \  
Q'_{-1} = 0, \ \ Q'_0 = b-1,   
$$
and 
$$
P'_{j+2} = \alpha'_{j+2} P_{j+1} + P'_j, \quad 
Q'_{j+2} = \alpha'_{j+2} Q_{j+1} + Q'_j, \quad j \ge -1.
$$
By construction, the sequence of pairs $((P'_j, Q'_j))_{j \ge 0}$ is a subsequence of $((P_j, Q_j))_{j \ge 0}$. 
Furthermore, it follows from (7.3) that $P'_j$ and $Q'_j$ are the numerator and 
denominator of 
\smallskip
$*$ the fraction $(1)_k$ if $\alpha'_j = c_k$; 

$*$ the fraction $(2)_k - (1)_k$ if $\alpha'_j = d_k$; 

$*$ the fraction $(2)_k$ if $\alpha'_j = 1$; 

$*$ the fraction $(3)_k$ if $\alpha'_j = e_k$; 

$*$ the fraction $(3)_{k+1}$ if $\alpha'_j = c_k + e_{k+1} + 1$; 

$*$ the fraction $(4)_k$ if $\alpha'_j = f_k$.
\smallskip

Now, we have to get rid of the $0$'s in $(\alpha'_j)_{j \ge 1}$ and construct 
a sequence $(\alpha''_j)_{j \ge 1}$ of positive integers. 
Since $e_k$ is positive for $k \ge 0$, there are no sequences of more than $3$ consecutive $0$'s 
in $(\alpha'_j)_{j \ge 1}$. 

As already observed, we have for nonnegative integers $x$ and $y$ we have
$$
\pmatrix{
x  & 1 \cr
1 & 0 \cr
}
\cdot
\pmatrix{
0 & 1 \cr
1 & 0 \cr
}
\cdot 
\pmatrix{
y & 1 \cr
1 & 0 \cr
}
=
\pmatrix{
x + y & 1 \cr
1 & 0 \cr
} 
$$
and, if $\alpha'_{j+1} = 0$, we get 
$$
\eqalign{
\pmatrix{
P'_{j+2} & P'_{j+1} \cr
Q'_{j+2} & Q'_{j +1} \cr
}
& =
\pmatrix{
P'_{j-1} & P'_{j-2} \cr
Q'_{j-1} & Q'_{j-2} \cr
}
\cdot
\pmatrix{
\alpha'_{j} & 1 \cr
1 & 0 \cr
} 
\cdot
\pmatrix{
\alpha'_{j+1} & 1 \cr
1 & 0 \cr
} 
\cdot
\pmatrix{
\alpha'_{j+2} & 1 \cr
1 & 0 \cr
} \cr
& = 
\pmatrix{
P'_{j-1} & P'_{j-2} \cr
Q'_{j-1} & Q'_{j-2} \cr
}
\cdot
\pmatrix{
\alpha'_{j} + \alpha'_{j+2} & 1 \cr
1 & 0 \cr
} \cr
& = 
\pmatrix{
(\alpha'_{j} + \alpha'_{j+2}) P'_{j-1} + P'_{j-2} & P'_{j-1} \cr
(\alpha'_{j} + \alpha'_{j+2}) Q'_{j - 1} + Q'_{j-2} & Q'_{j -1} \cr     
}. \cr
}   \eqno (7.4)
$$
This shows that $P'_{j-1}$ is followed by $P'_{j+2} = (\alpha'_{j} + \alpha'_{j+2}) P_{j-1} + P_{j-2}$, and   
similarly for $Q'_{j-1}$. 

By (7.4), if $x, 0, y$ are consecutive elements in $(\alpha'_j)_{j \ge 1}$, they have to be replaced 
by the single element $x+y$ in $(\alpha''_j)_{j \ge 1}$ and the pair associated with the 
partial quotients $x+y$ is the pair associated to $\alpha'_{j+2}$, that is, the pair $(P'_{j+2}, Q'_{j+2})$. 
Define recursively
$$
P''_{-1}= b-1, \  \ P''_0 = 0, \ \ 
Q''_{-1} = 0, \ \ Q''_0 = b-1, 
$$
and
$$
P''_{j+2} = \alpha''_{j+2} P''_{j+1} + P''_j, \quad 
Q''_{j+2} = \alpha''_{j+2} Q''_{j+1} + Q''_j, \quad j \ge -1.
$$
By construction, the sequence of pairs $((P''_j, Q''_j))_{j \ge 0}$ is a subsequence of $((P'_j, Q'_j))_{j \ge 0}$, 
hence of $((P_j, Q_j))_{j \ge 0}$. 
Let us discuss more in details which are the possible elements of the sequence $(\alpha''_j)_{j \ge 1}$. 
The following cases may occur:
\smallskip

$(i)$ $1, e_k, f_k=0, c_{k+1} + e_{k+2} + 1 \not=0, f_{k+2}$ are consecutive elements of $(\alpha'_j)_{j \ge 1}$, 
in which case we get the partial quotient $e_k + c_{k+1} + e_{k+2} + 1$ in $(\alpha''_j)_{j \ge 1}$ and
$$
[0; \alpha''_1, \ldots , e_k + c_{k+1} + e_{k+2} + 1] = (3)_{k+2},
$$
the preceding convergent being $(2)_k$.

$(ii)$ $1, e_k, f_k=0, c_{k+1} \not=0, d_{k+1}\not= 0$ are consecutive elements of $(\alpha'_j)_{j \ge 1}$, 
in which case we get the partial quotient $e_k + c_{k+1}$ in $(\alpha''_j)_{j \ge 1}$ and
$$
[0; \alpha''_1, \ldots , e_k + c_{k+1}] = (1)_{k+1},
$$
the preceding convergent being $(2)_k$.
 

$(iii)$ $e_k, f_k\not=0, c_{k+1} =0, d_{k+1} \not=0, 1$ are consecutive elements of $(\alpha'_j)_{j \ge 1}$, 
in which case we get the partial quotient  $f_k + d_{k+1}$ in $(\alpha''_j)_{j \ge 1}$ and
$$
[0; \alpha''_1, \ldots , e_k, f_k + d_{k+1}] = (2)_{k+1} - (1)_{k+1},
$$
the preceding convergent being $(3)_k$.


$(iv)$ $1, e_k, f_k=0, c_{k+1} \not=0, d_{k+1} = 0, 1$ are consecutive elements of $(\alpha'_j)_{j \ge 1}$, 
in which case we get the partial quotient $e_k + c_{k+1} + 1$  in $(\alpha''_j)_{j \ge 1}$ and
$$
[0; \alpha''_1, \ldots , e_k + c_{k+1} + 1] = (2)_{k+1},
$$
the preceding convergent being $(2)_k$.


$(v)$ $1, e_k, f_k \not=0, c_{k+1} \not=0, d_{k+1}= 0, 1$ are consecutive elements of $(\alpha'_j)_{j \ge 1}$, 
in which case we get the partial quotient $c_{k+1} + 1$ in $(\alpha''_j)_{j \ge 1}$ and
$$
[0; \alpha''_1, \ldots , c_{k+1} + 1] = (2)_{k+1},
$$
the preceding convergent being $(4)_k$.


$(vi)$ $1, e_k, f_k=0, c_{k+1} =0, d_{k+1} = 0, 1$ are consecutive elements of $(\alpha'_j)_{j \ge 1}$, 
in which case we get the partial quotient $e_k + 1$ in $(\alpha''_j)_{j \ge 1}$ and 
$$
[0; \alpha''_1, \ldots , e_k + 1] = (2)_{k+1},
$$
the preceding convergent being $(2)_k$.


\smallskip

The cases $(iv)$ to $(vi)$ occur only when $d_{k+1} = 0$, that is, when $b_1 = \ldots = b_{k+1} = 0$. 
They are not reflected in Proposition 6.2, where it is assumed that $t_k$ is positive.  

Note that if $f_k = c_{k+1} = 0$ or if $c_{k+1} = d_{k+1} = 0$, 
then there is nothing to do: we simply remove these two $0$'s from 
the sequence $(\alpha'_j)_{j \ge 1}$. 

\smallskip

The link with Proposition 6.2 is as follows:  

$*$ Case $(i)$ of Proposition 6.2 corresponds to the construction of $(\alpha'_j)_{j \ge 1}$ 
from $(\alpha_j)_{j \ge 1}$, with, if in addition $b_k$ is nul, Case $(i)$ above. 

$*$ Case $(ii)$ of Proposition 6.2 corresponds to Case $(iii)$ above if $b_{k+1}$ is positive, while 
if $b_{k+1}=0$ we remove two consecutive $0$'s in the sequence $(\alpha'_j)_{j \ge 1}$, thereby 
deleting two putative convergents.  

$*$ Case $(iii)$ of Proposition 6.2 corresponds to Case $(ii)$ above.   

\smallskip

Since the sequence $(\alpha''_j)_{j \ge 1}$ is composed of positive integers, the real number
$$
\zeta := [0; \alpha''_1, \alpha''_2, \ldots]
$$
is well defined by its continued fraction expansion. 
We have proved that all of its convergents are of the form $P_j / Q_j$ for some index $j$. 



It also follows from our discussion that $(2)_{k+1}$ is a convergent to $\zeta$ if $c_{k+1}$ 
and $c_{k+2}$ are nonnegative. If $c_{k+1} < 0$ and $f_{k-1} > 0$, then $f_{k-1}$ is an element of 
$(\alpha''_j)_{j \ge 1}$, associated with $(4)_{k-1} = (2)_{k+1}$. 
If $c_{k+1} < 0$ and $f_{k-1} = 0$, then $b_k = 0$ and
$1, e_{k-1} + c_{k} + e_{k+1} + 1$ are consecutive elements 
of $(\alpha''_j)_{j \ge 1}$, with this partial quotient $1$ being associated to $(2)_{k-1}$ and we have 
$(2)_{k-1} = (4)_{k-1} = (2)_{k+1}$. To summarize, we have shown that 
$(2)_{k+1}$ is a convergent to $\zeta$ unless $c_{k+2}$ is negative, that is, unless 
$a_{k+3} = b_{k+3}$. However, Proposition 6.1 asserts that 
$(2)_{k+1}$ is a convergent to $\xi$ if and only if
$a_{k+3} \ge b_{k+3} + 1$. Since there are infinitely many $h$ such that 
$a_{h} \ge b_{h} + 1$, we deduce that $\xi$ and $\zeta$ have infinitely many 
partial quotients in common, thus $\xi = \zeta$. 

The next statement summarizes what we have established.     
For $j \ge 1$, write $P_j/Q_j = [0;  \alpha''_1, \alpha''_2, \ldots, \alpha''_j]$ for the $j$-th 
convergent to $\xi$. 

\proclaim  
Proposition 7.2.  
All of the convergents to $\xi$ are of one of the five types
$(1)_{k}$, $(2)_{k} - (1)_k$, $(2)_k$, $(3)_k$, $(4)_k$. 
All its partial quotients are of the form 
$$
1, c_k, d_k, e_k, f_k, 
$$
or belong to the set 
$$
\{c_k + e_{k+1} + 1, e_{k-1} + c_k + e_{k+1} + 1, e_k + c_{k+1}, f_k + d_{k+1},
e_k + c_{k+1} + 1, e_k + 1, c_k + 1\}. 
$$
More precisely, we have
$$
P_j /Q_j = \cases{ 
(1)_k & if $\alpha''_j \in \{c_k, e_{k-1} + c_k\}$, \cr
(2)_k - (1)_k & if $\alpha''_j \in \{d_k, f_{k-1} + d_k \}$, \cr
(2)_k & if $\alpha''_j  \in \{1, e_{k-1} + c_k + 1, e_{k-1} + 1, c_k+1\}$,  \cr
(3)_k & if $\alpha''_j \in \{e_k, c_{k-1} + e_{k} + 1, e_{k-2} + c_{k-1} + e_{k} + 1\}$, \cr
(4)_k & if $\alpha''_j = f_k$. \cr}
$$

\vskip 5mm

\centerline{\bf 8. Remaining proofs}

\vskip 5mm

\noindent {\it Proof of Corollary 2.5. } 
Assume that $\theta$ has unbounded partial quotients 
(the case of bounded partial quotients is treated in Theorem 2.6). 
Let ${\cal K}$ be an infinite set of positive integers such that the subsequence 
$(a_k)_{k \in {\cal K}}$ is increasing. 
Assume first that there exists an infinite set $\cK' \subset \cK$ such that 
$(a_k - b_k)_{k \in \cK'}$ is increasing.  
For $k \ge 3$ in $\cK'$ we have
$$
\nu_{k-2} (4)= 1+ { r_{k}\over q_{k-1}} \ge  1+ { (a_k - b_k -1) q_{k-1} \over q_{k-1}},
$$
and, since $a_k - b_k$ can be arbitrarily large with $k$ in $\cK'$, we deduce that $\nu (4)$ is infinite. 


Assume now that there exist an infinite set $\cK' \subset \cK$    
and a nonnegative integer $\delta$ such that 
$a_k - b_k = \delta$ for $k$ in $\cK'$. For $k \ge 3$ in $\cK'$ we have
$$
\nu_{k-1} (3) = 1+ { q_{k}\over r_{k}+ q_{k-1}} \ge { a_k q_{k-1} \over r_{k-1}+ \delta q_{k-1} + q_{k-2}}
\ge {a_k \over \delta + 2}.
$$
We deduce that $\nu(3)$ is infinite. Consequently, any Sturmian number whose slope has 
unbounded partial quotients is a Liouville number. 
\cqfd

\noindent {\it Proof of Theorem 2.6. } 
Assume that $\theta$ has bounded partial quotients. 
Observe that 
$$
\nu_k(3) = 1+ { q_{k+1}\over r_{k+1}+ q_k} \le 1+ { q_{k+1}\over  q_k},
\quad \nu_k(4)= 1+ { r_{k+2}\over q_{k+1}} \le 1+ { q_{k+2}\over q_{k+1}}.
$$
If $a_{k+1} = b_{k+1}$, then $r_{k+1} = r_{k-1}$ and $t_k = t_{k-1}$, thus 
$$
\nu_k(2) = 2+ {r_{k} \over r_{k+1}+ t_k} = 2 + {r_{k} \over r_{k-1}+ t_{k-1}} \le 2 + {q_{k} \over q_{k-1}}.
$$
If $a_{k+1} > b_{k+1}$, then $r_{k+1} \ge r_k + q_{k-1}$, thus
$$
\nu_k(2) = 2+ {r_{k} \over r_{k+1}+ t_k} \le 2 + {r_{k} \over r_{k}+ t_{k}} \le 3, 
$$
and
$$
\nu_k(1) = 2+ {t_k\over r_{k+1}} \le 2+ {t_k \over r_k + q_{k-1}}  \le 2+ {q_k \over q_{k-1}}.
$$
This shows that the irrationality exponent of $\xi_b (\theta, \rho)$ satisfies
$$
\mu (\xi_b (\theta, \rho)) \le 2 + \limsup_{k \to + \infty} \, {q_k \over q_{k-1}} = 1 + \mu (\xi_b (\theta)). 
\eqno (8.1) 
$$

Let us now show that there exist intercepts $\rho$ for which equality holds. 
Let $\cK$ be an infinite set of positive integers such that 
$$
\lim_{k \to + \infty, k \in \cK} \, {q_k \over q_{k-1}} =   \mu (\xi_b (\theta)) - 1. 
$$
Take $k_1 \ge 3$ in $\cK$ and
set $b_1 = \ldots = b_{k_1} = 0$. Put $a_{k_1 + 1} = b_{k_1 + 1}$ and 
$b_{k_1 + 2} = b_{k_1 + 3} = \ldots = b_{k_2} = 0$, 
where $k_2 > k_1 + 2$ is in $\cK$ and sufficiently large to ensure that $r_{k_2} \ge q_{k_2} / 2$. 
Then, put $b_{k_2 + 1} = a_{k_2 + 1}$ and $b_{k_2 + 2} = \ldots = b_{k_3} = 0$,
where $k_3 > k_2 + 2$ is in $\cK$ and sufficiently large to ensure that $r_{k_3} \ge 2 q_{k_3} / 3$. 
Proceeding like this, we define inductively 
an icreasing sequence $(k_j)_{j \ge 2}$ of integers in $\cK$ such that $b_{k_j + 1} = a_{k_j + 1}$
and $b_k = 0$ for every $k$ not in $(k_j)_{j \ge 2}$. In addition, we have 
$r_{k_j} \ge (j-1) q_{k_j} / j$, for $j \ge 2$.  

Denote by $\rho$ the intercept defined by this sequence $(b_k)_{k \ge 1}$ and let us 
determine the irrationality exponent of $\xi_b (\theta, \rho)$. 

Recall that for an index $h$ such that $b_{h+1} = a_{h+1}$ we have $r_{h+1} = r_{h-1}$ 
and $t_h = t_{h-1}$, thus
$$
\nu_h (2) = 2+ {r_{h} \over r_{h+1}+ t_h} = 2 + {r_{h} \over r_{h-1}+ t_{h-1}} 
= 2 + {r_{h} \over q_{h-1}}. 
$$
Consequently, we get
$$
\nu_{k_j} (2) \ge 2 + { (j-1) q_{k_j} \over j q_{k_j-1}}, \quad j \ge 2,
$$
and 
$$
\mu (\xi_b (\theta, \rho))  \ge \nu (2) \ge 2 + \limsup_{j \to + \infty} \, {q_{k_j} \over q_{k_j-1}} 
= 2 + \lim_{k \to + \infty, k \in \cK} \, {q_k \over q_{k-1}} =   \mu (\xi_b (\theta)) + 1. 
$$
The reverse inequality follows from (8.1). 
Consequently, we get
$$
\mu (\xi_b (\theta, \rho))  = 1 + \mu (\xi_b (\theta)). 
$$
This proves the theorem. 
\cqfd

\noindent {\it Proof of Theorem 2.7. } 
Assume that not all $b_k$ are $0$. 
Let $k$ be an integer large enough to ensure that $t_k$ is positive and that 
$a_k, a_{k+1}, \ldots $ are all at most equal to $M$. 
Then, it follows from Proposition 6.2 that there are four (possibly overlapping) cases:
\smallskip
$(i)$ If $a_{k+2} = b_{k+2}$, then $(2)_{k+1}$ is a convergent to $\xi$;

$(ii)$ If $a_{k+3} = b_{k+3}$, then $(2)_{k+2}$ is a convergent to $\xi$;

$(iii)$ If $a_{k+4} = b_{k+4}$, then $(2)_{k+3}$ is a convergent to $\xi$;

$(iv)$ If $(i)$, $(ii)$, and $(iii)$ do not hold, then $(1)_{k+1}$ and $(2)_{k+1}$ are convergents to $\xi$.
\smallskip

In case $(i)$, the rate of approximation of $\xi$ by $(2)_{k+1}$ is at least equal to 
$$
\nu_{k+1} (2) = 2+ {r_{k+1} \over r_{k+2}+ t_{k+1}} 
= 2 + {r_{k+1} \over q_{k}} \ge 2 + {q_{k-1} \over q_{k}} \ge 2 + {1 \over M + 1},
$$
since $r_{k+2} = r_k$, $t_{k+1} = t_k$, and $r_{k+1} \ge q_{k-1}$. 

Similarly, in case $(ii)$ (resp., $(iii)$), the rate of approximation of $\xi$ by $(2)_{k+2}$ 
(resp., by $(2)_{k+3}$) is at least equal to $2 + 1/(M+1)$. 

In case $(iv)$, note that $r_{k+2} + t_{k+1} \le (a_{k+2} + 1) q_{k+1}$, thus
$$
\nu_{k+1} (1) = 2+ {t_{k+1} \over r_{k+2}} \ge 2+ {t_{k+1} \over (M + 1) q_{k+1}}, \quad 
\nu_{k+1} (2) = 2+ {r_{k+1} \over r_{k+2}+ t_{k+1}} \ge 2+ {r_{k+1} \over (M + 1) q_{k+1}}. 
$$
Recalling that $r_{k+1} + t_{k+1} = q_{k+1}$, we get 
$$
\max\{\nu_{k+1} (1), \nu_{k+1} (2)\} \ge 2+ {1 \over 2 (M + 1)}. 
$$
This shows that, for every sufficiently large $k$, there exists a rational number $P / Q$ with 
$$
b^{q_{k}} \le Q \le b^{q_{k+4}}    \eqno (8.2)
$$
such that $|\xi - P /Q | \le Q^{-2 - 1/ (2 (M+1))}$. 
We are then in position to apply Th\'eor\`eme 3.1 of \cite{AdBu10} with $\eps = {1 \over 2 (M + 1)}$
and ${\cal S}$ the empty set. Note that, by (8.2), the number $c$ introduced in $(3.2)$ of its proof 
can be taken to be $(M+1)^5$. Consequently, the upper bound 
$$
w_d^* (\xi) \le (2d)^{\kappa (\log \log 3d)}, \quad d \ge 1, 
$$
given by Th\'eor\`eme 3.1 of \cite{AdBu10} holds with a real number $\kappa$ 
depending only on $M$.

\vskip 8mm

\centerline{\bf References}

\vskip 5mm

\beginthebibliography{999}

\bibitem{Ad10}
B. Adamczewski, 
{\it On the expansion of some exponential periods in an integer base},
Math. Ann. 346 (2010), 107--116. 

\bibitem{AdAl07}
B. Adamczewski and J.-P. Allouche,
{\it Reversals and palindromes in continued fractions},
Theor. Comput. Sci.  380  (2007),  220--237.

\bibitem{AdBu10}
B. Adamczewski et Y. Bugeaud,
{\it Mesures de transcendance et aspects quantitatifs de la 
m\'ethode de Thue--Siegel--Roth--Schmidt},
Proc. London Math. Soc. 101 (2010), 1--31.

\bibitem{AdBu11}
B. Adamczewski and Y. Bugeaud,
{\it Nombres r\'eels de complexit\'e sous-lin\'eaire : 
mesures d'irrationalit\'e et de transcendance},
J. Reine Angew. Math. 658 (2011), 65--98.

\bibitem{AlSh03}
J.-P. Allouche and J. Shallit, 
Automatic Sequences: Theory, Applications, Generalizations,  
Cambridge University Press, 2003.

\bibitem{AdDa77}
W. W. Adams and J. L. Davison,
{\it A remarkable class of continued fractions},
Proc. Amer. Math. Soc. 65 (1977), 194--198.

\bibitem{Ar02}
P. Arnoux,  
{\it Sturmian sequences}. 
In: Substitutions in dynamics, arithmetics and combinatorics, 143--198,
Lecture Notes in Math., 1794, Springer, Berlin, 2002. 

\bibitem{BHZ06}
V. Berth\'e, C. Holton, and L. Q. Zamboni,
{\it Initial powers of Sturmian sequences},
Acta Arith.  122  (2006),  315--347.

\bibitem{Ber01}
V. Berth\'e,
{\it Autour du syst\`eme de num\'eration d'Ostrowski}, 
Bull. Belg. Math. Soc. Simon Stevin 8 (2001),  209--239.

\bibitem{Bohm27}
P. E. B\"ohmer,
{\it \"Uber die Transzendenz gewisser dyadischer Br\"uche},
Math. Ann. 96 (1927), 367--377. 

\bibitem{BuLiv1}
Y. Bugeaud,
Approximation by algebraic numbers. 
Cambridge Tracts in Mathematics 160,
Cambridge, 2004.

\bibitem{BuLiv2}
Y. Bugeaud,
Distribution modulo one and Diophantine approximation.
Cambridge Tracts in Mathematics 193, Cambridge, 2012. 

\bibitem{BuKim19}
Y. Bugeaud and D. H. Kim, 
{\it A new complexity function, repetitions in Sturmian words, 
and irrationality exponents of Sturmian numbers},
Trans. Amer. Math. Soc. 371 (2019), 3281--3308. 

\bibitem{BKLN}
Y. Bugeaud, D. H. Kim, M. Laurent and A. Nogueira, 
{\it On the Diophantine nature of the elements of Cantor sets arising in the dynamics of contracted rotations}, 
Ann. Scuola Normale Superiore di Pisa.
To appear. 
{\tt https://arxiv.org/abs/2001.00380}   

\bibitem{Dan72}
L. V. Danilov,
{\it Certain classes of transcendental numbers}
Mat. Zametki 12 (1972), 149--154 (in Russian). 
English translation in Math. Notes 12 (1972), 524--527.

\bibitem{Davi77}
J. L. Davison,
{\it A series and its associated continued fraction},
Proc. Amer. Math. Soc. 63 (1977), 29--32.

\bibitem{FeMa97}
S. Ferenczi and Ch. Mauduit,
{\it Transcendence of numbers with a low complexity expansion},
J. Number Theory 67 (1997), 146--161.

\bibitem{Ko39}
J. F. Koksma,
{\it \"Uber die Mahlersche Klasseneinteilung der transzendenten Zahlen
und die Approximation komplexer Zahlen durch algebraische Zahlen},
Monats. Math. Phys. 48 (1939), 176--189.

\bibitem{Kom96}
T. Komatsu, 
{\it A certain power series and the inhomogeneous continued fraction expansions}, 
J. Number Theory 59 (1996), 291--312. 

\bibitem{LaNo18}
M. Laurent and A. Nogueira,
{\it Rotation number of contracted rotations}, 
J. Mod. Dyn. 12 (2018), 175--191.

\bibitem{LaNo21}
M. Laurent and A. Nogueira,
{\it Dynamics of 2-interval piecewise affine maps and Hecke-Mahler series}, 
J. Mod. Dyn. 17 (2021), 33--63.

\bibitem{Loth02}
M. Lothaire, 
Algebraic combinatorics on words. 
Encyclopedia of Mathematics and its Applications, vol. 90, Cambridge University Press, Cambridge, 2002.

\bibitem{Wo19}
C. Wojcik,
Factorisations des mots de basse complexit\'e. 
Doctoral thesis, Universit\'e de Lyon, 2019.

\endthebibliography

\vskip1cm

\noindent Yann Bugeaud      \hfill    Michel Laurent

\noindent Universit\'e de Strasbourg, CNRS    \hfill  Aix-Marseille Universit\'e, CNRS, Centrale Marseille 

\noindent IRMA, UMR 7501     \hfill   Institut de Math\'ematiques de Marseille

\noindent 7, rue Ren\'e Descartes         \hfill    163 avenue de Luminy, Case 907

\noindent 67084 STRASBOURG  (FRANCE)   \hfill  13288  MARSEILLE C\'edex 9 (FRANCE)

\vskip2mm

\noindent {\tt bugeaud@math.unistra.fr}    \hfill {\tt michel-julien.laurent@univ-amu.fr} 

\vskip1cm

\bye

\noindent Michel Laurent

\noindent Aix-Marseille Universit\'e, CNRS, Centrale Marseille

\noindent Institut de Math\'ematiques de Marseille

\noindent 163 avenue de Luminy, Case 907

\noindent 13288  MARSEILLE C\'edex 9 (FRANCE)

\vskip2mm

\noindent{\tt michel-julien.laurent@univ-amu.fr}

\bye